# STOCHASTIC INTEGRATION IN UMD BANACH SPACES

By J. M. A. M. van Neerven ,[1,2] M. C. Veraar [1] and L. Weis [3]

*Delft University of Technology, Delft University of Technology and Technische Universität Karlsruhe*

In this paper we construct a theory of stochastic integration of processes with values in $\mathcal{L}(H,E)$, where $H$ is a separable Hilbert space and $E$ is a UMD Banach space (i.e., a space in which martingale differences are unconditional). The integrator is an $H$-cylindrical Brownian motion. Our approach is based on a two-sided $L^p$-decoupling inequality for UMD spaces due to Garling, which is combined with the theory of stochastic integration of $\mathcal{L}(H,E)$-valued functions introduced recently by two of the authors. We obtain various characterizations of the stochastic integral and prove versions of the Itô isometry, the Burkholder–Davis–Gundy inequalities, and the representation theorem for Brownian martingales.

**1. Introduction.** It is well known that the theory of stochastic integration can be extended to Hilbert space-valued processes in a very satisfactory way. The reason for this is that the Itô isometry is an $L^2$-isometry which easily extends to the Hilbert space setting. At the same time, this explains why it is considerably more difficult to formulate a theory of stochastic integration for processes taking values in a Banach space $E$. By a well-known result due to Rosiński and Suchanecki [36], the class of strongly measurable functions $\phi:[0,T] \to E$ that are stochastically integrable (in a sense that is made precise below) with respect to a Brownian motion $W$ coincides with $L^2(0,T;E)$ if and only if $E$ isomorphic to a Hilbert space. More precisely,

Received April 2005; revised October 2006.
[1]Supported in part by a "VIDI subsidie" (639.032.201) in the "Vernieuwingsimpuls" program of the Netherlands Organization for Scientific Research (NWO).
[2]Supported in part by the Research Training Network "Evolution Equations for Deterministic and Stochastic Systems" (HPRN-CT-2002-00281).
[3]Supported by grants from the Volkswagenstiftung (I/78593) and the Deutsche Forschungsgemeinschaft (We 2847/1-2).

the authors showed that $E$ has type 2 if and only if every $\phi \in L^2(0,T;E)$ is stochastically integrable and there is a constant $C \geq 0$ such that

$$\mathbb{E}\left\|\int_0^T \phi(t)\, dW(t)\right\|^2 \leq C^2 \|\phi\|_{L^2(0,T;E)}^2,$$

and that $E$ has cotype 2 if and only if every strongly measurable, stochastically integrable function $\phi$ belongs to $L^2(0,T;E)$ and there exists a constant $C \geq 0$ such that

$$\|\phi\|_{L^2(0,T;E)}^2 \leq C^2 \mathbb{E}\left\|\int_0^T \phi(t)\, dW(t)\right\|^2.$$

Combined with Kwapień's theorem which asserts that $E$ is isomorphic to a Hilbert space if and only if $E$ has both type 2 and cotype 2, this gives the result as stated.

It turns out that the Itô isometry does extend to the Banach space setting provided one reformulates it properly. To this end let us first observe that, for Hilbert spaces $E$,

$$\|\phi\|_{L^2(0,T;E)} = \|I_\phi\|_{\mathcal{L}_2(L^2(0,T),E)},$$

where $\mathcal{L}_2(L^2(0,T),E)$ denotes the space of Hilbert–Schmidt operators from $L^2(0,T)$ to $E$ and $I_\phi \colon L^2(0,T) \to E$ is the integral operator defined by

$$I_\phi f := \int_0^T f(t)\phi(t)\, dt.$$

Now one observes that the class $\mathcal{L}_2(L^2(0,T),E)$ coincides isometrically with the class of $\gamma$-radonifying operators $\gamma(L^2(0,T),E)$. With this in mind one has the natural result that a function $\phi \colon [0,T] \to E$, where $E$ is now an arbitrary Banach space, is stochastically integrable if and only if the corresponding integral operator $I_\phi$ defines an element in $\gamma(L^2(0,T),E)$, and if this is the case the Itô isometry takes the form

$$\mathbb{E}\left\|\int_0^T \phi(t)\, dW(t)\right\|^2 = \|I_\phi\|_{\gamma(L^2(0,T),E)}^2.$$

This operator-theoretic approach to stochastic integration of $E$-valued functions has been developed systematically by two of us [28]. The purpose of the present paper is to extend this theory to the case of $E$-valued processes. This is achieved by the decoupling approach initiated by Garling [15], who proved a two-sided $L^p$-estimate for the stochastic integral of an elementary adapted process $\phi$ with values in a UMD space in terms of the stochastic integral of $\phi$ with respect to an independent Brownian motion. A new short proof of these estimates is included. The decoupled integral is defined path by path, which makes it possible to apply the theory developed for $E$-valued functions to the sample paths of $\phi$. As a result, we obtain a two-sided $L^p$-estimate for the stochastic integral of $\phi$ in terms of the $L^p$-norm of the associated $\gamma(L^2(0,T),E)$-valued random variable $X_\phi$ defined path by path by $X_\phi(\omega) := I_{\phi(\cdot,\omega)}$. As it turns out, the space $L^p(\Omega;\gamma(L^2(0,T),E))$ provides the right setting to establish a fairly complete theory of stochastic integration of adapted processes with values in a UMD space $E$. We obtain various characterizations of the class of stochastically integrable processes and prove a version of the Itô isometry, which, together with Doob's maximal inequality, leads to the following Burkholder–Davis–Gundy type inequalities: for every $p \in (1,\infty)$ there exist constants $0 < c < C < \infty$, depending only on $p$ and $E$, such that

$$(1.1) \quad c^p \mathbb{E}\|X_\phi\|_{\gamma(L^2(0,T),E)}^p \leq \mathbb{E} \sup_{t \in [0,T]} \left\| \int_0^t \phi(s)\, dW(s) \right\|^p \leq C^p \mathbb{E}\|X_\phi\|_{\gamma(L^2(0,T),E)}^p.$$

This result clearly indicates that for UMD spaces $E$, the space $L^p(\Omega;\gamma(L^2(0,T),E))$ is the correct space of integration, at least if one is interested in having two-sided $L^p$-estimates for the stochastic integral. In order to keep this paper at a reasonable length, the proof of an Itô formula is postponed to the paper [26].

The fact that the two-sided estimates (1.1) are indeed available shows that our theory extends the Hilbert space theory in a very natural way. Garling's estimates actually characterize the class of UMD spaces, and for this reason the decoupling approach is naturally limited to this class of spaces if one insists on having two-sided estimates. From the point of view of applications this is an acceptable limitation, since this class includes many of the classical reflexive spaces such as the $L^p$ spaces for $p \in (1,\infty)$ as well as spaces constructed from these, such as Sobolev spaces and Besov spaces. At the price of obtaining only one-sided estimates, our theory can be extended to a class of Banach spaces having a one-sided randomized version of the UMD property. This class of spaces was introduced by Garling [16] and includes all $L^1$-spaces. The details will be presented elsewhere.

For the important special case of $L^q(S)$-spaces, where $(S,\Sigma,\mu)$ is a $\sigma$-finite measure space and $q \in (1,\infty)$, the operator language can be avoided and the norm of $L^p(\Omega;\gamma(L^2(0,t),L^q(S)))$ is equivalent to a square function norm. More precisely, for every $p \in (1,\infty)$ there exist constants $0 < c < C < \infty$ such that

$$c^p \mathbb{E} \left\| \left( \int_0^T |\phi(t,\cdot)|^2\, dt \right)^{1/2} \right\|_{L^q(S)}^p \leq \mathbb{E}\|X_\phi\|_{\gamma(L^2(0,T),E)}^p \leq C^p \mathbb{E} \left\| \left( \int_0^T |\phi(t,\cdot)|^2\, dt \right)^{1/2} \right\|_{L^q(S)}^p.$$

As an application of our abstract results we prove in Section 4 that $L^p$-martingales with values in a UMD space are stochastically integrable and we provide an estimate for their stochastic integrals.

A decoupling inequality for the moments of tangent martingale difference sequences was obtained by Hitczenko [17] and McConnell [24]. McConnell used it to obtain sufficient pathwise conditions for stochastic integrability of processes with values in a UMD space. We shall recover McConnell's result by localization. This approach has the advantage of replacing the $\zeta$-convexity arguments used by McConnell by abstract operator-theoretic arguments. In our approach, the UMD property is only used through the application of Garling's estimates which are derived directly from the definition of the UMD property.

With only little extra effort the results described above can be derived in the more general setting of $\mathcal{L}(H,E)$-valued processes, with $H$-cylindrical Brownian motions as integrators. Here, $H$ is a separable real Hilbert space and $\mathcal{L}(H,E)$ denotes the space of bounded linear operators from $H$ to $E$. We shall formulate all results in this framework, because this permits the application of our theory to the study of certain classes of nonlinear stochastic evolution equations in $E$, driven by an $H$-cylindrical Brownian motion. Here the space $L^p(\Omega; \gamma(L^2(0,T;H),E))$ [which takes over the role of $L^p(\Omega; \gamma(L^2(0,T),E))$] serves as the framework for a classical fixed point argument. This will be the topic of a forthcoming paper [27]. The reader who is not interested in this level of generality may simply substitute $H$ by $\mathbb{R}$ and identify $\mathcal{L}(\mathbb{R},E)$ with $E$ and $W_H$ with a Brownian motion $W$ throughout the paper.

Many authors (cf. [1, 4, 5, 6, 11, 12, 30, 31] and references therein) have considered the problem of stochastic integration in Banach spaces with martingale type 2 or related geometric properties. We compare their approaches with ours at the end of Section 3. Various classical spaces, such as $L^q(S)$ for $q \in (1,2)$, do have the UMD property but fail to have martingale type 2. On the other hand, an example due to Bourgain [2] implies the existence of martingale type 2 spaces which do not have the UMD property.

Preliminary versions of this paper have been presented at the meeting Stochastic Partial Differential Equations and Applications-VII in Levico Terme in January 2004 (M. V.) and meeting Spectral Theory in Banach Spaces and Harmonic Analysis in Oberwolfach in July 2004 (J. v. N.).

**2. Operator-valued processes.** Throughout this paper, $(\Omega, \mathcal{F}, \mathbb{P})$ is a probability space endowed with a filtration $\mathbb{F} = (\mathcal{F}_t)_{t \in [0,T]}$ satisfying the usual conditions, $H$ is a separable real Hilbert space, and $E$ is a real Banach space with dual $E^*$. The inner product of two elements $h_1, h_2 \in H$ is written as $[h_1, h_2]_H$, and the duality pairing of elements $x \in E$ and $x^* \in E^*$ is denoted by $\langle x, x^* \rangle$. We use the notation $\mathcal{L}(H,E)$ for the space of all bounded linear operators from $H$ to $E$. We shall always identify $H$ with its dual in the natural way. In particular, the adjoint of an operator in $\mathcal{L}(H,E)$ is an operator in $\mathcal{L}(E^*,H)$.

We write $Q_1 \lesssim_A Q_2$ to express that there exists a constant $c$, only depending on $A$, such that $Q_1 \leq cQ_2$. We write $Q_1 \eqsim_A Q_2$ to express that $Q_1 \lesssim_A Q_2$ and $Q_2 \lesssim_A Q_1$.

2.1. *Measurability.* Let $(S, \Sigma)$ be a measurable space and let $E$ be a real Banach space with dual space $E^*$. A function $f: S \to E$ is called *measurable* if $f^{-1}(B) \in \Sigma$ for every Borel set $B \subseteq E$, and *simple* if it is measurable and takes finitely many values. The function $f$ is called *strongly measurable* if it is the pointwise limit of a sequence of simple functions, and *separably valued* if there exists a separable closed subspace $E_0$ of $E$ such that $f(s) \in E_0$ for all $s \in S$. Given a functional $x^* \in E^*$, we define the function $\langle f, x^* \rangle: S \to \mathbb{R}$ by $\langle f, x^* \rangle(s) := \langle f(s), x^* \rangle$. The function $f$ is said to be *scalarly measurable* if $\langle f, x^* \rangle$ is measurable for all $x^* \in E^*$. More generally, if $F$ is a linear subspace of $E^*$ and $\langle f, x^* \rangle$ is measurable for all $x^* \in F$, we say that $f$ is $F$-*scalarly measurable*. The following result is known as the Pettis measurability theorem ([37], Proposition I.1.10).

PROPOSITION 2.1 (Pettis measurability theorem). *For a function $f: S \to E$ the following assertions are equivalent:*

(1) *$f$ is strongly measurable;*
(2) *$f$ is separably valued and scalarly measurable;*
(3) *$f$ is separably valued and $F$-scalarly measurable for some weak$^*$-dense linear subspace $F$ of $E^*$.*

A function $\Phi: S \to \mathcal{L}(H, E)$ is called *scalarly measurable* if the function $\Phi^* x^*: S \to H$ defined by $\Phi^* x^*(s) := \Phi^*(s) x^*$ is strongly measurable for all $x^* \in E^*$, and $H$-*strongly measurable* if for all $h \in H$ the function $\Phi h: S \to E$ defined by $\Phi h(s) := \Phi(s)h$ is strongly measurable.

Let $\mu$ be a finite measure on $(S, \Sigma)$. Two scalarly measurable functions $\Phi, \Psi: S \to \mathcal{L}(H, E)$ are called *scalarly $\mu$-equivalent* if for all $x^* \in E^*$ we have $\Phi^* x^* = \Psi^* x^*$ $\mu$-almost everywhere on $S$.

PROPOSITION 2.2. *If $E$ is weakly compactly generated, then every scalarly measurable function $\Phi: S \to \mathcal{L}(H, E)$ is scalarly $\mu$-equivalent to an $H$-strongly measurable function $\Psi: S \to \mathcal{L}(H, E)$.*

For $H = \mathbb{R}$ this is a deep result of [14], and the result for general $H$ is easily deduced from it. Recall that a Banach space $E$ is *weakly compactly* *generated* if it is the closed linear span of one of its weakly compact subsets. All separable Banach spaces and all reflexive Banach spaces are weakly compactly generated.

In the main results of this paper we are concerned with $\mathcal{L}(H, E)$-valued stochastic processes $(\Phi_t)_{t \in [0,T]}$ on a probability space $(\Omega, \mathcal{F}, \mathbb{P})$, which will be viewed as functions $\Phi : [0, T] \times \Omega \to \mathcal{L}(H, E)$. Since $E$ will always be a Banach space belonging to a certain class of reflexive Banach spaces, Proposition 2.2 justifies us to restrict our considerations to $H$-strongly measurable processes, that is, to processes $\Phi : [0, T] \times \Omega \to \mathcal{L}(H, E)$ with the property that for all $h \in H$ the $E$-valued process $\Phi h : [0, T] \times \Omega \to E$ defined by $\Phi h(t, \omega) := \Phi(t, \omega)h$ is strongly measurable. We point out, however, that most of our proofs work equally well for scalarly measurable processes.

2.2. $\gamma$-*Radonifying operators.* In this subsection we discuss some properties of the operator ideal of $\gamma$-radonifying operators from a separable real Hilbert space $\mathcal{H}$ to $E$. The special case $\mathcal{H} = L^2(0, T; H)$ will play an important role in this paper.

Let $(\gamma_n)_{n \geq 1}$ be a sequence of independent standard Gaussian random variables on a probability space $(\Omega', \mathcal{F}', \mathbb{P}')$ [we reserve the notation $(\Omega, \mathcal{F}, \mathbb{P})$ for the probability space on which our processes live] and let $\mathcal{H}$ be a separable real Hilbert space. A bounded operator $R \in \mathcal{L}(\mathcal{H}, E)$ is said to be $\gamma$-*radonifying* if there exists an orthonormal basis $(h_n)_{n \geq 1}$ of $\mathcal{H}$ such that the Gaussian series $\sum_{n \geq 1} \gamma_n R h_n$ converges in $L^2(\Omega'; E)$. We then define

$$\|R\|_{\gamma(\mathcal{H}, E)} := \left( \mathbb{E}' \left\| \sum_{n \geq 1} \gamma_n R h_n \right\|^2 \right)^{1/2}.$$

This number does not depend on the sequence $(\gamma_n)_{n \geq 1}$ and the basis $(h_n)_{n \geq 1}$, and it defines a norm on the space $\gamma(\mathcal{H}, E)$ of all $\gamma$-radonifying operators from $\mathcal{H}$ into $E$. Endowed with this norm, $\gamma(\mathcal{H}, E)$ is a Banach space, which is separable if $E$ is separable. If $R \in \gamma(\mathcal{H}, E)$, then $\|R\| \leq \|R\|_{\gamma(\mathcal{H}, E)}$. If $E$ is a Hilbert space, then $\gamma(\mathcal{H}, E) = \mathcal{L}_2(\mathcal{H}, E)$ isometrically, where $\mathcal{L}_2(\mathcal{H}, E)$ denotes the space of all Hilbert–Schmidt operators from $H$ to $E$.

The following property of $\gamma$-radonifying operators will be important:

PROPOSITION 2.3 (Ideal property). *Let $\tilde{E}$ be a real Banach space and let $\tilde{\mathcal{H}}$ be a separable real Hilbert space. If $B_1 \in \mathcal{L}(\tilde{\mathcal{H}}, \mathcal{H})$, $R \in \gamma(\mathcal{H}, E)$ and $B_2 \in \mathcal{L}(E, \tilde{E})$, then $B_2 \circ R \circ B_1 \in \gamma(\tilde{\mathcal{H}}, \tilde{E})$ and $\|B_2 \circ R \circ B_1\|_{\gamma(\tilde{\mathcal{H}}, \tilde{E})} \leq \|B_2\| \|R\|_{\gamma(\mathcal{H}, E)} \|B_1\|$.*

For these and related results we refer to [13, 30, 37].
We shall frequently use the following convergence result.

PROPOSITION 2.4. *If the $T_1, T_2, \ldots \in \mathcal{L}(\mathcal{H})$ and $T \in \mathcal{L}(\mathcal{H})$ satisfy:*

(1) $\sup_{n \geq 1} \|T_n\| < \infty$,
(2) $T^* h = \lim_{n \to \infty} T_n^* h$ *for all* $h \in \mathcal{H}$,

*then for all $R \in \gamma(\mathcal{H}, E)$ we have $R \circ T = \lim_{n \to \infty} R \circ T_n$ in $\gamma(\mathcal{H}, E)$.*

PROOF. By the estimate $\|R \circ S\|_{\gamma(\mathcal{H}, E)} \leq \|R\|_{\gamma(\mathcal{H}, E)} \|S\|$ for $S \in \mathcal{L}(H)$ and (1), it suffices to consider finite rank operators $R \in \gamma(\mathcal{H}, E)$. For such an operator, say $R = \sum_{j=1}^{k} h_j \otimes x_j$, we may estimate

$$\|R \circ (T - T_n)\|_{\gamma(\mathcal{H}, E)} \leq \sum_{j=1}^{k} \|x_j\| \|T^* h_j - T_n^* h_j\|.$$

By (2), the right-hand side tends to zero as $n \to \infty$. □

Identifying $\mathcal{H} \otimes E^*$ canonically with a weak$^*$-dense linear subspace of $(\gamma(\mathcal{H}, E))^*$, as an easy consequence of the Pettis measurability theorem we obtain the following measurability result for $\gamma(\mathcal{H}, E)$-valued functions. A closely related result is given in [30].

LEMMA 2.5. *Let $(S, \Sigma, \mu)$ be a $\sigma$-finite measure space. For a function $X : S \to \gamma(\mathcal{H}, E)$ the following assertions are equivalent:*

(1) *The function $s \mapsto X(s)$ is strongly measurable;*
(2) *For all $h \in \mathcal{H}$, the function $s \mapsto X(s)h$ is strongly measurable.*

*If these equivalent conditions hold, there exists a separable closed subspace $E_0$ of $E$ such that $X(s) \in \gamma(\mathcal{H}, E_0)$ for all $s \in S$.*

The following result will be useful:

PROPOSITION 2.6 ($\gamma$-Fubini isomorphism). *Let $(S, \Sigma, \mu)$ be a $\sigma$-finite measure space and let $p \in [1, \infty)$ be fixed. The mapping $F_\gamma : L^p(S; \gamma(\mathcal{H}, E)) \to \mathcal{L}(\mathcal{H}, L^p(S; E))$ defined by*

$$(F_\gamma(X)h)(s) := X(s)h, \qquad s \in S, h \in \mathcal{H},$$

*defines an isomorphism from $L^p(S; \gamma(\mathcal{H}, E))$ onto $\gamma(\mathcal{H}, L^p(S; E))$.*

PROOF. Let $(h_n)_{n \geq 1}$ be an orthonormal basis for $\mathcal{H}$ and let $(\gamma_n)_{n \geq 1}$ be a sequence of independent standard Gaussian random variables on a probability space $(\Omega', \mathcal{F}', \mathbb{P}')$. By the Kahane–Khinchine inequalities and Fubini's theorem we have, for any $X \in L^p(S; \gamma(\mathcal{H}, E))$,

$$\|F_\gamma(X)\|_{\gamma(\mathcal{H}, L^p(S;E))}$$

$$= \left(\mathbb{E}'\left\|\sum_{n\geq 1}\gamma_n F_\gamma(X)h_n\right\|^2_{L^p(S;E)}\right)^{1/2} \eqsim_p \left(\mathbb{E}'\left\|\sum_{n\geq 1}\gamma_n F_\gamma(X)h_n\right\|^p_{L^p(S;E)}\right)^{1/p}$$

(2.1)

$$= \left(\int_S \mathbb{E}'\left\|\sum_{n\geq 1}\gamma_n X h_n\right\|^p d\mu\right)^{1/p} \eqsim_p \left(\int_S \left(\mathbb{E}'\left\|\sum_{n\geq 1}\gamma_n X h_n\right\|^2\right)^{p/2} d\mu\right)^{1/p}$$

$$= \left(\int_S \|X\|^p_{\gamma(\mathcal{H},E)} d\mu\right)^{1/p} = \|X\|_{L^p(S;\gamma(\mathcal{H},E))}.$$

By these estimates the range of the operator $X \mapsto F_\gamma(X)$ is closed in $\gamma(\mathcal{H}, L^p(S; E))$. Hence to show that this operator is surjective it is enough to show that its range is dense. But this follows from

$$F_\gamma\left(\sum_{n=1}^N \mathbf{1}_{S_n} \otimes \left(\sum_{k=1}^K h_k \otimes x_{kn}\right)\right) = \sum_{k=1}^K h_k \otimes \left(\sum_{n=1}^N \mathbf{1}_{S_n} \otimes x_{kn}\right),$$

for all $S_n \in \Sigma$ with $\mu(S_n) < \infty$ and $x_{kn} \in E$, noting that the elements on the right-hand side are dense in $\gamma(\mathcal{H}, E)$. □

For $p = 2$ we have equality in all steps of (2.1).

For later use we note that if $(S, \Sigma, \mu) = (\Omega, \mathcal{F}, \mathbb{P})$ is a probability space and $\mathcal{H} = L^2(0, T; H)$, then the $\gamma$-Fubini isomorphism takes the form

$$F_\gamma : L^p(\Omega; \gamma(L^2(0, T; H), E)) \simeq \gamma(L^2(0, T; H), L^p(\Omega; E)).$$

The space on the left-hand side will play an important role in the stochastic integration theory developed in Section 3.

2.3. *Representation.* As before we let $H$ is a separable real Hilbert space.

An $H$-strongly measurable function $\Phi : [0, T] \to \mathcal{L}(H, E)$ is said to *belong to $L^2(0, T; H)$ scalarly* if for all $x^* \in E^*$ the function $\Phi^* x^* : (0, T) \to H$ belongs to $L^2(0, T; H)$. Such a function *represents* an operator $R \in \mathcal{L}(L^2(0, T; H), E)$ if for all $f \in L^2(0, T; H)$ and $x^* \in E^*$ we have

$$\langle Rf, x^* \rangle = \int_0^T \langle \Phi(t)f(t), x^* \rangle\, dt.$$

Similarly, an $H$-strongly measurable process $\Phi : [0, T] \times \Omega \to \mathcal{L}(H, E)$ is said to *belong to $L^2(0, T; H)$ scalarly almost surely* if for all $x^* \in E^*$ it is true that the function $\Phi^*_\omega x^* : (0, T) \to E$ belongs to $L^2(0, T; H)$ for almost all $\omega \in \Omega$. Here we use the notation

$$\Phi_\omega(t) := \Phi(t, \omega).$$

Note that the exceptional set may depend on $x^*$. Such a process $\Phi$ is said to *represent* an $H$-strongly measurable random variable $X: \Omega \to \mathcal{L}(L^2(0,T;H), E)$ if for all $f \in L^2(0,T;H)$ and $x^* \in E^*$ we have

$$\langle X(\omega)f, x^* \rangle = \int_0^T [f(t), \Phi^*(t,\omega)x^*]_H \, dt \qquad \text{for almost all } \omega \in \Omega.$$

If $\Phi_1$ and $\Phi_2$ are $H$-strongly measurable, then $\Phi_1$ and $\Phi_2$ represent the same random variable $X$ if and only if $\Phi_1(t,\omega) = \Phi_2(t,\omega)$ for almost all $(t,\omega) \in [0,T] \times \Omega$. In the converse direction, the strongly measurable random variables $X_1$ and $X_2$ are represented by the same process $\Phi$ if and only if $X_1(\omega) = X_2(\omega)$ for almost all $\omega \in \Omega$.

For a random variable $X: \Omega \to \gamma(L^2(0,T;H), E)$ we denote by $\langle X, x^* \rangle : \Omega \to L^2(0,T;H)$ the random variable defined by

$$\langle X, x^* \rangle(\omega) := X^*(\omega)x^*.$$

Notice that $X$ is represented by $\Phi$ if and only if for all $x^* \in E^*$, $\langle X, x^* \rangle = \Phi^* x^*$ in $L^2(0,T;H)$ almost surely.

The next lemma relates the above representability concepts and shows that the exceptional sets may be chosen independently of $x^*$.

LEMMA 2.7. *Let $\Phi: [0,T] \times \Omega \to \mathcal{L}(H,E)$ be an $H$-strongly measurable process and let $X: \Omega \to \gamma(L^2(0,T;H), E)$ be strongly measurable. The following assertions are equivalent:*

(1) *$\Phi$ represents $X$.*
(2) *$\Phi_\omega$ represents $X(\omega)$ for almost all $\omega \in \Omega$.*

PROOF. The implication (1) $\Rightarrow$ (2) is clear from the definitions. To prove the implication (2) $\Rightarrow$ (1) we start by noting that the Pettis measurability theorem allows us to assume, without loss of generality, that $E$ is separable. Let $(f_m)_{m \geq 1}$ be a dense sequence in $L^2(0,T;H)$ and let $(x_n^*)_{n \geq 1}$ be a sequence in $E^*$ with weak$^*$-dense linear span. Choose a null set $N \subseteq \Omega$ such that:

(i) $\Phi^*(\cdot, \omega)x_n^* \in L^2(0,T;H)$ for all $x_n^*$ and all $\omega \in \complement N$;
(ii) for all $f_m$, all $x_n^*$, and all $\omega \in \complement N$,

(2.2) $$\langle X(\omega)f, x^* \rangle = \int_0^T \langle \Phi(t,\omega)f(t), x^* \rangle \, dt.$$

Let $F$ denote the linear subspace of all $x^* \in E^*$ for which:

(i)$'$ $\Phi^*(\cdot, \omega)x^* \in L^2(0,T;H)$ for all $\omega \in \complement N$;
(ii)$'$ (2.2) holds for all $f \in L^2(0,T;H)$ and all $\omega \in \complement N$.

By a limiting argument we see that $x_n^* \in F$ for all $n \geq 1$. Hence $F$ is weak*-dense. We claim that $F$ is also weak*-sequentially closed. Once we have checked this, we obtain $F = E^*$ by the Krein–Smulyan theorem, see [7], Proposition 1.2.

To prove the claim, fix $\omega \in \complement N$ and $x^* \in F$ arbitrary. Then, by (2.2),

$$(2.3) \qquad \|\Phi^*(\cdot,\omega)x^*\|_{L^2(0,T;H)} \leq \|X(\omega)\|_{\gamma(L^2(0,T;H),E)}\|x^*\|.$$

Suppose now that $\lim_{n\to\infty} y_n^* = y^*$ weak* in $E^*$ with each $y_n^* \in F$. Then (2.3) shows that the sequence $\Phi^*(\cdot,\omega)y_n^*$ is bounded in $L^2(0,T;H)$. By a convex combination argument as in [7], Proposition 2.2, we find that $y^* \in F$, and the claim is proved. □

REMARK 2.8. The assumptions of (2) already imply that the induced mapping $\omega \mapsto X(\omega)$ from $\Omega$ to $\gamma(L^2(0,T;H),E)$ has a strongly measurable version. To see this, first note that by Lemma 2.5 it suffices to show that for all $f \in L^2(0,T;H)$ the mapping $\omega \mapsto X(\omega)f$ is strongly measurable from $\Omega$ to $E$. By assumption, almost surely we have that (2.2) holds for all $f \in L^2(0,T;H)$ and $x^* \in E^*$. By the $H$-strong measurability of $\Phi$ and Fubini's theorem, the right-hand side of (2.2) is a measurable function of $\omega$. Thus $\omega \mapsto X(\omega)f$ is scalarly measurable. By the Pettis measurability theorem it remains to show that $\omega \mapsto X(\omega)f$ is almost surely separably-valued.

Since $t \mapsto \Phi(t,\omega)$ is $H$-strongly measurable for almost all $\omega \in \Omega$ and belongs to $L^2(0,T;H)$ scalarly, it follows that $t \mapsto \Phi(t,\omega)f(t)$ is Pettis integrable with

$$X(\omega)f = \int_0^T \Phi(t,\omega)f(t)\,dt$$

for almost all $\omega \in \Omega$. Then by the Hahn–Banach theorem, $\omega \mapsto X(\omega)f$ almost surely takes its values in the closed subspace spanned by the range of $(t,\omega) \mapsto \Phi(t,\omega)f(t)$, which is separable by the $H$-strong measurability of $\Phi$.

The following example shows what might go wrong if the assumption of representation in Lemma 2.7 were to be replaced by the weaker assumption of belonging to $L^2(0,T;H)$ scalarly almost surely, even in the simple case where $H = \mathbb{R}$ and $E$ is a separable real Hilbert space.

EXAMPLE 2.9. Let $E$ be an infinite-dimensional separable Hilbert space with inner product $[\cdot,\cdot]_E$. We shall construct a process $\phi \colon [0,1] \times \Omega \to E$ with the following properties:

(1) $\phi$ is strongly measurable;
(2) $\phi$ belongs to $L^2(0,1)$ scalarly almost surely;
(3) $\phi_\omega$ fails to be scalarly in $L^2(0,1)$ for almost all $\omega \in \Omega$.

Let $(\xi_n)_{n\geq 1}$ denote a sequence of independent $\{0,1\}$-valued random variables on a probability space $(\Omega, \mathcal{F}, \mathbb{P})$ satisfying $\mathbb{P}\{\xi_n = 1\} = \frac{1}{n}$ for $n \geq 1$. Fix an orthonormal basis $(x_n)_{n\geq 1}$ in $E$. Define $\phi\colon [0,1] \times \Omega \to E$ by $\phi(0,\omega) = 0$ and

$$\phi(t,\omega) := n^{1/2} 2^{n/2} \xi_n(\omega) x_n \qquad \text{for } n \geq 1 \text{ and } t \in [2^{-n}, 2^{-n+1}).$$

It is clear that $\phi$ is strongly measurable, and (2) is checked by direct computation. To check (3) we first note that

$$\mathbb{P}\{\xi_n = 1 \text{ for infinitely many } n \geq 1\} = 1.$$

Indeed, this follows from the fact that for each $n \geq 1$ we have

$$\mathbb{P}\{\xi_k = 0 \text{ for all } k \geq n\} = \prod_{k \geq n}\left(1 - \frac{1}{k}\right) = 0.$$

Fix an arbitrary $\omega \in \Omega$ for which $\xi_n = 1$ for infinitely many $n \geq 1$, say $\xi_n(\omega) = 1$ for $n = n_1, n_2, \ldots$ and $\xi_n(\omega) = 0$ otherwise. Let $(a_n)_{n\geq 1}$ be any sequence of real numbers with $\sum_{n\geq 1} a_n^2 < \infty$ and $\sum_{n\geq 1} n a_n^2 = \infty$, and put $x := \sum_{k\geq 1} a_k x_{n_k}$. Then,

$$\int_0^1 [\phi(t,\omega), x]_E^2 \, dt = \sum_{k \geq 1} n_k a_k^2 \geq \sum_{k \geq 1} k a_k^2 = \infty.$$

This concludes the construction.

### 2.4. Adaptedness.

A process $\Phi\colon [0,T] \times \Omega \to \mathcal{L}(H,E)$ is said to be *elementary adapted* to the filtration $\mathbb{F} = (\mathcal{F}_t)_{t \in [0,T]}$ if it is of the form

$$(2.4) \qquad \Phi(t,\omega) = \sum_{n=0}^{N} \sum_{m=1}^{M} \mathbf{1}_{(t_{n-1}, t_n] \times A_{mn}}(t,\omega) \sum_{k=1}^{K} h_k \otimes x_{kmn},$$

where $0 \leq t_0 < \cdots < t_N \leq T$ and the sets $A_{1n}, \ldots, A_{Mn} \in \mathcal{F}_{t_{n-1}}$ are disjoint for each $n$ (with the understanding that $(t_{-1}, t_0] := \{0\}$ and $\mathcal{F}_{t_{-1}} := \mathcal{F}_0$) and the vectors $h_1, \ldots, h_K \in H$ are orthonormal. An $H$-strongly measurable process $\Phi\colon [0,T] \times \Omega \to \mathcal{L}(H,E)$ is called *adapted* to $\mathbb{F}$ if for all $h \in H$ the $E$-valued process $\Phi h$ is *strongly adapted*, that is, for all $t \in [0,T]$ the random variable $\Phi(t)h$ is strongly $\mathcal{F}_t$-measurable.

A random variable $X\colon \Omega \to \gamma(L^2(0,T;H), E)$ is *elementary adapted* to $\mathbb{F}$ if it is represented by an elementary adapted process. We call $X$ *strongly adapted* to $\mathbb{F}$ if there exists a sequence of elementary adapted random variables $X_n\colon \Omega \to \gamma(L^2(0,T;H), E)$ such that $\lim_{n\to\infty} X_n = X$ in measure in $\gamma(L^2(0,T;H), E)$.

Recall that for a finite measure space $(S, \Sigma, \mu)$ and strongly measurable functions $f, f_1, f_2, \ldots$ from $S$ into a Banach space $F$, $f = \lim_{n \to \infty} f_n$ in measure if and only if $\lim_{n\to\infty} \mathbb{E}(\|f - f_n\|_B \wedge 1) = 0$.

PROPOSITION 2.10. *For a strongly measurable random variable $X:\Omega \to \gamma(L^2(0,T;H),E)$, the following assertions are equivalent:*

(1) *$X$ is strongly adapted to $\mathbb{F}$;*
(2) *$X(\mathbf{1}_{[0,t]}f)$ is strongly $\mathcal{F}_t$-measurable for all $f \in L^2(0,T;H)$ and $t \in [0,T]$.*

PROOF. The implication $(1) \Rightarrow (2)$ follows readily from the definitions.

$(2) \Rightarrow (1)$: For $\delta \geq 0$ we define the right translate $R^\delta$ of an operator $R \in \gamma(L^2(0,T;H),E)$ by

$$R^\delta f := Rf_\delta, \qquad f \in L^2(0,T;H),$$

where $f_\delta$ denotes the left translate of $f$. It follows by the right ideal property and Proposition 2.4 that $R^\delta \in \gamma(L^2(0,T;H),E)$ with $\|R^\delta\|_{\gamma(H,E)} \leq \|R\|_{\gamma(H,E)}$ and that $\delta \mapsto R^\delta$ is continuous with respect to the $\gamma$-radonifying norm.

Define the right translate $X^\delta: \Omega \to \gamma(L^2(0,T;H),E)$ by pointwise action, that is, $X^\delta(\omega) := (X(\omega))^\delta$. Note that $X^\delta$ is strongly measurable by Lemma 2.5. By dominated convergence, $\lim_{\delta \downarrow 0} X^\delta = X$ in measure in $\gamma(L^2(0,T;H),E)$. Thus, for $\varepsilon > 0$ fixed, we may choose $\delta > 0$ such that

(2.5) $$\mathbb{E}(\|X - X^\delta\|_{\gamma(L^2(0,T;H),E)} \wedge 1) < \varepsilon.$$

Let $0 = t_0 < \cdots < t_N = T$ be an arbitrary partition of $[0,T]$ of mesh $\leq \delta$ and let $I_n = (t_{n-1}, t_n]$ for $n = 1, \ldots, N$. Let $X_n^\delta$ denote the restriction of $X^\delta$ to $I_n$, that is,

$$X_n^\delta(\omega)g := X^\delta(\omega)i_n g, \qquad g \in L^2(I_n; H),$$

where $i_n: L^2(I_n; H) \to L^2(0,T;H)$ is the inclusion mapping. From the assumption (1) we obtain that $X_n^\delta$ is strongly $\mathcal{F}_{t_{n-1}}$-measurable as a random variable with values in $\gamma(L^2(I_n;H),E))$. Pick a simple $\mathcal{F}_{t_{n-1}}$-measurable random variable $Y_n: \Omega \to \gamma(L^2(I_n;H),E)$ such that

$$\mathbb{E}(\|X_n^\delta - Y_n\|_{\gamma(L^2(I_n;H),E)} \wedge 1) < \frac{\varepsilon}{N},$$

say $Y_n = \sum_{m=1}^{M_n} \mathbf{1}_{A_{mn}} \otimes S_{mn}$ with $A_{mn} \in \mathcal{F}_{t_{n-1}}$ and $S_{mn} \in \gamma(L^2(I_n;H),E)$. By a further approximation we may assume that the $S_{mn}$ are represented by elementary functions $\Psi_{mn}: [0,T] \to \mathcal{L}(H,E)$ of the form

$$\Psi_{mn}(t) = \sum_{j=1}^{J_{mn}} \mathbf{1}_{(s_{(j-1)mn}, s_{jmn}]}(t) \sum_{k=1}^{K_{mn}} (h_k \otimes x_{kmn}),$$

where $t_{n-1} \leq s_{0mn} < \cdots < s_{J_{mn}mn} \leq t_n$ and $(h_k)_{k\geq 1}$ is a fixed orthonormal basis for $H$. Define the process $\Psi: [0,T] \times \Omega \to \mathcal{L}(H,E)$ by

$$\Psi(t,\omega) := \sum_{m=1}^{M_n} \mathbf{1}_{A_{mn}}(\omega) \Psi_{mn}(t), \qquad t \in I_n.$$

It is easily checked that $\Psi$ is elementary adapted. Let $Y : \Omega \to \gamma(L^2(0, T; H), E)$ be represented by $\Psi$. Then $Y$ is elementary adapted and satisfies

$$(2.6) \qquad \mathbb{E}(\|X^\delta - Y\|_{\gamma(L^2(0,T;H),E)} \wedge 1) < \varepsilon.$$

Finally, by (2.5) and (2.6),

$$\mathbb{E}(\|X - Y\|_{\gamma(L^2(0,T;H),E)} \wedge 1) \leq 2\varepsilon.$$

This proves that $X$ can be approximated in measure by a sequence of elementary adapted elements $X_n$. □

PROPOSITION 2.11. *If $\Phi : [0, T] \times \Omega \to \mathcal{L}(H, E)$ is an $H$-strongly measurable and adapted process representing a random variable $X : \Omega \to \gamma(L^2(0, T; H), E)$, then $X$ is strongly adapted to $\mathbb{F}$.*

PROOF. By using the identity $\langle X(\mathbf{1}_{[0,t]} f), x^* \rangle = [\mathbf{1}_{[0,t]} f, \Phi^* x^*]_{L^2(0,T;H)}$ and noting that the right-hand side is $\mathcal{F}_t$-measurable, this follows trivially from Proposition 2.10 and the Pettis measurability theorem. □

For $p \in [1, \infty)$, the closure in $L^p(\Omega; \gamma(L^2(0, T; H), E))$ of the elementary adapted elements will be denoted by

$$L^p_\mathbb{F}(\Omega; \gamma(L^2(0, T; H), E)).$$

PROPOSITION 2.12. *If the random variable $X \in L^p(\Omega; \gamma(L^2(0, T; H), E))$ is strongly adapted to $\mathbb{F}$, then $X \in L^p_\mathbb{F}(\Omega; \gamma(L^2(0, T; H), E))$.*

PROOF. By assumption, condition (1) in Proposition 2.10 is satisfied. Now we can repeat the proof of the implication $(1) \Rightarrow (2)$, but instead of approximating in measure we approximate in the $L^p$-norm. □

**3. $L^p$-stochastic integration.** Recall that a family $W_H = (W_H(t))_{t \in [0, T]}$ of bounded linear operators from $H$ to $L^2(\Omega)$ is called an *$H$-cylindrical Brownian motion* if:

(1) $W_H h = (W_H(t) h)_{t \in [0,T]}$ is real-valued Brownian motion for each $h \in H$,

(2) $\mathbb{E}(W_H(s) g \cdot W_H(t) h) = (s \wedge t)[g, h]_H$ for all $s, t \in [0, T], g, h \in H$.

We always assume that the $H$-cylindrical Brownian motion $W_H$ is adapted to a given filtration $\mathbb{F}$ satisfying the usual conditions, that is, the Brownian motions $W_H h$ are adapted to $\mathbb{F}$ for all $h \in H$.

EXAMPLE 3.1.   Let $W = (W(t))_{t \geq 0}$ be an $E$-valued Brownian motion and let $C \in \mathcal{L}(E^*, E)$ be its covariance operator, that is, $C$ is the unique positive symmetric operator such that $\mathbb{E}\langle W(t), x^* \rangle^2 = t\langle Cx^*, x^* \rangle$ for all $t \geq 0$ and $x^* \in E^*$. Let $H_C$ be the reproducing kernel Hilbert space associated with $C$ and let $i_C : H_C \hookrightarrow E$ be the inclusion operator. Then the mappings
$$W_{H_C}(t) : i_C^* x^* \mapsto \langle W(t), x^* \rangle$$
uniquely extend to an $H_C$-cylindrical Brownian motion $W_{H_C}$.

If $\Phi : [0, T] \times \Omega \to E$ is an elementary adapted process of the form (2.4), we define the stochastic integral $\int_0^T \Phi(t) \, dW_H(t)$ by
$$\int_0^T \Phi(t) \, dW_H(t) := \sum_{n=1}^N \sum_{m=1}^M \mathbf{1}_{A_{mn}} \sum_{k=1}^K (W_H(t_n) h_k - W_H(t_{n-1}) h_k) x_{kmn}.$$
Note that the stochastic integral belongs to $L^p(\Omega; E)$ for all $p \in [1, \infty)$. It turns out that for a suitable class of Banach spaces $E$ this definition can be extended to the class of adapted processes representing an element of $L^p(\Omega; \gamma(L^2(0, T; H), E))$. In order to motivate our approach, we recall the following result on stochastic integration of $\mathcal{L}(H, E)$-valued functions from [28]; see [7, 23, 35, 36] for related results.

PROPOSITION 3.2.    *For a function $\Phi : [0, T] \to \mathcal{L}(H, E)$ belonging to $L^2(0, T; H)$ scalarly, the following assertions are equivalent:*

  (1)  *There exists a sequence $(\Phi_n)_{n \geq 1}$ of elementary functions such that:*

  (i)  *for all $x^* \in E^*$ we have $\lim_{n \to \infty} \Phi_n^* x^* = \Phi^* x^*$ in $L^2(0, T; H)$,*
  (ii) *there exists a strongly measurable random variable $\eta : \Omega \to E$ such that*
$$\eta = \lim_{n \to \infty} \int_0^T \Phi_n(t) \, dW_H(t) \qquad \text{in probability};$$

  (2)  *There exists a strongly measurable random variable $\eta : \Omega \to E$ such that for all $x^* \in E^*$ we have*
$$\langle \eta, x^* \rangle = \int_0^T \Phi^*(t) x^* \, dW_H(t) \qquad \text{almost surely};$$

  (3)  *$\Phi$ represents an operator $R \in \gamma(L^2(0, T; H), E)$.*

*In this situation the random variables $\eta$ in (1) and (2) are uniquely determined and equal almost surely. Moreover, $\eta$ is Gaussian and for all $p \in [1, \infty)$ we have*

(3.1) $\qquad (\mathbb{E}\|\eta\|^p)^{1/p} \eqsim_p (\mathbb{E}\|\eta\|^2)^{1/2} = \|R\|_{\gamma(L^2(0,T;H), E)}.$

*For all $p \in [1, \infty)$ the convergence in (1), part (ii), is in $L^p(\Omega; E)$.*

A function $\Phi$ satisfying the equivalent conditions of Proposition 3.2 will be called *stochastically integrable* with respect to $W_H$. The random variable $\eta$ is called the *stochastic integral* of $\Phi$ with respect to $W_H$, notation

$$\eta =: \int_0^T \Phi(t)\, dW_H(t).$$

The second identity in (3.1) may be interpreted as an analogue of the Itô isometry.

REMARK 3.3. If $\Phi$ is $H$-strongly measurable and belongs to $L^2(0,T;H)$ scalarly, the arguments in [36] can be adapted to show that condition (1) is equivalent to

(1)′ There exists a sequence $(\Phi_n)_{n\geq 1}$ of elementary functions such that:

(i) for all $h \in H$ we have $\lim_{n\to\infty} \Phi_n h = \Phi h$ in measure on $[0,T]$,

(ii) there exists a strongly measurable random variable $\eta: \Omega \to E$ such that

$$\eta = \lim_{n\to\infty} \int_0^T \Phi_n(t)\, dW_H(t) \qquad \text{in probability.}$$

The extension of Proposition 3.2 to processes is based on a decoupling inequality for processes with values in a UMD space $E$. Recall that a Banach space $E$ is a *UMD space* if for some (equivalently, for all) $p \in (1,\infty)$ there exists a constant $\beta_{p,E} \geq 1$ such that for every $n \geq 1$, every martingale difference sequence $(d_j)_{j=1}^n$ in $L^p(\Omega;E)$, and every $\{-1,1\}$-valued sequence $(\varepsilon_j)_{j=1}^n$ we have

$$\left(\mathbb{E}\left\|\sum_{j=1}^n \varepsilon_j d_j\right\|^p\right)^{1/p} \leq \beta_{p,E} \left(\mathbb{E}\left\|\sum_{j=1}^n d_j\right\|^p\right)^{1/p}.$$

Examples of UMD spaces are all Hilbert spaces and the spaces $L^p(S)$ for $1 < p < \infty$ and $\sigma$-finite measure spaces $(S, \Sigma, \mu)$. If $E$ is a UMD space, then $L^p(S;E)$ is a UMD space for $1 < p < \infty$. For an overview of the theory of UMD spaces we refer the reader to [8, 34] and references given therein.

Let $\tilde{W}_H$ be an $H$-cylindrical Brownian motion on a second probability space $(\tilde{\Omega}, \tilde{\mathcal{F}}, \tilde{\mathbb{P}})$, adapted to a filtration $\tilde{\mathbb{F}}$. If $\Phi: [0,T] \times \Omega \to E$ is an elementary adapted process of the form (2.4), we define the *decoupled* stochastic integral $\int_0^T \Phi(t)\, d\tilde{W}_H(t)$ by

$$\int_0^T \Phi(t)\, d\tilde{W}_H(t) := \sum_{n=1}^N \sum_{m=1}^M \mathbf{1}_{A_{mn}} \sum_{k=1}^K (\tilde{W}_H(t_n)h_k - \tilde{W}_H(t_{n-1})h_k) x_{kmn}.$$

This stochastic integral belongs $L^p(\Omega; L^p(\tilde{\Omega}; E))$.

The following result was proved by Garling [15], Theorems 2 and 2′, for finite-dimensional Hilbert spaces $H$. For reasons of completeness we include a short proof which is a variation of a more general argument in [25].

LEMMA 3.4 (Decoupling). *Let $H$ be a nonzero separable real Hilbert space and fix $p \in (1, \infty)$. The following assertions are equivalent:*

(1) *$E$ is a UMD space;*
(2) *For every elementary adapted process $\Phi : [0, T] \times \Omega \to \mathcal{L}(H, E)$ we have*

$$\beta_{p,E}^{-p} \mathbb{E}\tilde{\mathbb{E}} \left\| \int_0^T \Phi(t)\, d\tilde{W}_H(t) \right\|^p \leq \mathbb{E} \left\| \int_0^T \Phi(t)\, dW_H(t) \right\|^p$$
$$\leq \beta_{p,E}^p \mathbb{E}\tilde{\mathbb{E}} \left\| \int_0^T \Phi(t)\, d\tilde{W}_H(t) \right\|^p.$$

PROOF. $(1) \Rightarrow (2)$: Let $\Phi$ be an elementary adapted process of the form (2.4). We extend $\Phi$, as well as $W_H$, $\tilde{W}_H$ and the $\sigma$-algebras $\mathcal{F}_t$, $\tilde{\mathcal{F}}_t$ in the obvious way to $\Omega \times \tilde{\Omega}$. Write

$$\sum_{n=1}^N d_n = \int_0^T \Phi(t)\, dW_H(t) \quad \text{and} \quad \sum_{n=1}^N e_n = \int_0^T \Phi(t)\, d\tilde{W}_H(t),$$

where the random variables $d_n$ and $e_n$ on $\Omega \times \tilde{\Omega}$ are defined by $d_n = W_H(t_n)\xi_n - W_H(t_{n-1})\xi_n$ and $e_n = \tilde{W}_H(t_n)\xi_n - \tilde{W}_H(t_{n-1})\xi_n$, where $\xi_n := \sum_{m=1}^M \mathbf{1}_{A_{mn}} \times \sum_{k=1}^K h_k \otimes x_{kmn}$ and

$$W_H(t)\xi_n := \sum_{m=1}^M \mathbf{1}_{A_{mn}} \sum_{k=1}^K W_H(t)h_k \otimes x_{kmn}.$$

For $n = 1, \ldots, N$ let

$$r_{2n-1} := \tfrac{1}{2}(d_n + e_n) \quad \text{and} \quad r_{2n} := \tfrac{1}{2}(d_n - e_n).$$

Then, $(r_j)_{j=1}^{2N}$ is a martingale difference sequence with respect to the filtration $(\mathcal{G}_j)_{j=1}^{2N}$, where

$$\mathcal{G}_{2n} = \sigma(\mathcal{F}_{t_n} \otimes \tilde{\mathcal{F}}_{t_n}),$$
$$\mathcal{G}_{2n-1} = \sigma(\mathcal{F}_{t_{n-1}} \otimes \tilde{\mathcal{F}}_{t_{n-1}}, w_{n1}, w_{n2}, \ldots),$$

where

$$w_{nk} = (W_H(t_n)h_k - W_H(t_{n-1})h_k) + (\tilde{W}_H(t_n)h_k - \tilde{W}_H(t_{n-1})h_k).$$

Notice that

$$\sum_{n=1}^N d_n = \sum_{j=1}^{2N} r_j \quad \text{and} \quad \sum_{n=1}^N e_n = \sum_{j=1}^{2N} (-1)^{j+1} r_j.$$

Hence (2) follows from the UMD property applied to the sequences $(r_j)_{j=1}^{2N}$ and $((-1)^{j+1} r_j)_{j=1}^{2N}$.

(2) $\Rightarrow$ (1): See [15], Theorem 2. $\square$

If $X \in L^p(\Omega; \gamma(L^2(0,T;H), E))$ is elementary adapted, we define the random variable $I^{W_H}(X) \in L^p(\Omega; E)$ by
$$I^{W_H}(X) := \int_0^T \Phi(t) \, dW_H(t),$$
where $\Phi$ is an elementary adapted process representing $X$. Note that $I^{W_H}(X)$ does not depend on the choice of the representing process $\Phi$. Clearly $I^{W_H}(X) \in L_0^p(\Omega, \mathcal{F}_T; E)$, the closed subspace of $L^p(\Omega; E)$ consisting of all $\mathcal{F}_T$-measurable random variables with mean zero. In the first main result of this section we extend the mapping $X \mapsto I^{W_H}(X)$ to a bounded operator from $L_{\mathbb{F}}^p(\Omega; \gamma(L^2(0,T;H), E))$ to $L_0^p(\Omega, \mathcal{F}_T; E)$. If $\mathbb{F} = \mathbb{F}^{W_H}$ is the augmented filtration generated by the Brownian motions $W_H h$, $h \in H$, this mapping turns out to be an isomorphism.

THEOREM 3.5 (Itô isomorphism). *Let $E$ be a UMD space and fix $p \in (1, \infty)$. The mapping $X \mapsto I^{W_H}(X)$ has a unique extension to a bounded operator*
$$I^{W_H} : L_{\mathbb{F}}^p(\Omega; \gamma(L^2(0,T;H), E)) \to L_0^p(\Omega, \mathcal{F}_T; E).$$
*This operator is an isomorphism onto its range and we have the two-sided estimate*
$$\beta_{p,E}^{-p} \|X\|_{L^p(\Omega; \gamma(L^2(0,T;H),E))} \lesssim_p \mathbb{E}\|I^{W_H}(X)\|^p \lesssim_p \beta_{p,E} \|X\|_{L^p(\Omega; \gamma(L^2(0,T;H),E))}^p.$$
*For the augmented Brownian filtration $\mathbb{F}^{W_H}$ we have an isomorphism of Banach spaces*
$$I^{W_H} : L_{\mathbb{F}^{W_H}}^p(\Omega; \gamma(L^2(0,T;H), E)) \simeq L_0^p(\Omega, \mathcal{F}_T^{W_H}; E).$$

PROOF. Let $X \in L^p(\Omega; \gamma(L^2(0,T;H), E))$ be elementary and adapted, and let $\Phi$ be an elementary adapted process representing $X$. It follows from Proposition 3.2, the Kahane–Khinchine inequalities and Lemma 3.4 that
$$\mathbb{E}\|X\|_{\gamma(L^2(0,T;H),E)}^p = \mathbb{E}\left\|\int_0^T \Phi(t) \, d\tilde{W}_H(t)\right\|_{L^2(\tilde{\Omega}; E)}^p$$
$$\eqsim_p \mathbb{E}\left\|\int_0^T \Phi(t) \, d\tilde{W}_H(t)\right\|_{L^p(\tilde{\Omega}; E)}^p$$
$$\eqsim_{p,E} \mathbb{E}\left\|\int_0^T \Phi(t) \, dW_H(t)\right\|^p = \mathbb{E}\|I^{W_H}(X)\|^p.$$

Thus the map $X \mapsto I^{W_H}(X)$ extends uniquely to an isomorphism from $L^p_{\mathbb{F}}(\Omega; \gamma(L^2(0,T;H),E))$ onto its range, which is a closed subspace of $L^p_0(\Omega, \mathcal{F}_0; E)$.

Next assume that $\mathbb{F} = \mathbb{F}^{W_H}$. Since $I^{W_H}$ is an isomorphism onto its range, which is a closed subspace of $L^p_0(\Omega, \mathcal{F}^{W_H}_T; E)$, it suffices to show that this operator has dense range in $L^p_0(\Omega, \mathcal{F}^{W_H}_T; E)$.

Let $(h_k)_{k \geq 1}$ be a fixed orthonormal basis for $H$. For $m = 1, 2, \ldots$ let $\mathcal{F}^{(m)}_T$ be denote by the augmented $\sigma$-algebra generated by $\{W_H(t)h_k : t \in [0,T], 1 \leq k \leq m\}$. Since $\mathcal{F}^{W_H}_T$ is generated by the $\sigma$-algebras $\mathcal{F}^{(m)}_T$, by the martingale convergence theorem and approximation we may assume $\eta$ is in $L^p_0(\Omega, \mathcal{F}^{(m)}_T; E)$ and of the form $\sum_{n=1}^N (\mathbf{1}_{A_n} - P(A_n)) \otimes x_n$ with $A_n \in \mathcal{F}^m_T$ and $x_n \in E$. From linearity and the identity

$$I^{W_H}(\phi \otimes x) = (I^{W_H}(\phi)) \otimes x, \qquad \phi \in L^p_{\mathbb{F}}(\Omega; L^2(0,T;H)),$$

it even suffices to show that $\mathbf{1}_{A_n} - P(A_n) = I^{W_H}(\phi)$ for some $\phi \in L^p_{\mathbb{F}}(\Omega; L^2(0,T;H))$. By the Itô representation theorem for Brownian martingales (cf. [18], Lemma 18.11 and [20], Theorem 3.4.15), there exists $\phi \in L^2_{\mathbb{F}}(\Omega; L^2(0,T;H))$ such that $\mathbf{1}_{A_n} - P(A_n) = \int_0^T \phi(t)\,dW(t)$, and the Burkholder–Davis–Gundy inequalities and Doob's maximal inequality imply that $\phi \in L^p_{\mathbb{F}}(\Omega; L^2(0,T;H))$. □

We return to the general setting where $W_H$ is adapted to an arbitrary filtration $\mathbb{F}$ satisfying the usual conditions. The second main result of this section describes the precise relationship between the $L^p$-stochastic integral and the operator $I^{W_H}$. It extends Proposition 3.2 to $\mathcal{L}(H,E)$-valued processes. In view of Proposition 2.2 we restrict ourselves to $H$-strongly measurable processes.

THEOREM 3.6.  *Let $E$ be a UMD space and fix $p \in (1, \infty)$. For an $H$-strongly measurable and adapted process $\Phi : [0,T] \times \Omega \to \mathcal{L}(H,E)$ belonging to $L^p(\Omega; L^2(0,T;H))$ scalarly, the following assertions are equivalent:*

(1) *There exists a sequence $(\Phi_n)_{n \geq 1}$ of elementary adapted processes such that:*

(i) *for all $h \in H$ and $x^* \in E^*$ we have $\lim_{n \to \infty} \langle \Phi_n h, x^* \rangle = \langle \Phi h, x^* \rangle$ in measure on $[0,T] \times \Omega$,*

(ii) *there exists a strongly measurable random variable $\eta \in L^p(\Omega; E)$ such that*

$$\eta = \lim_{n \to \infty} \int_0^T \Phi_n(t)\,dW_H(t) \qquad in\ L^p(\Omega; E);$$

(2) *There exists a strongly measurable random variable* $\eta \in L^p(\Omega; E)$ *such that for all* $x^* \in E^*$ *we have*

$$\langle \eta, x^* \rangle = \int_0^T \Phi^*(t) x^* \, dW_H(t) \qquad \text{in } L^p(\Omega);$$

(3) $\Phi$ *represents an element* $X \in L^p(\Omega; \gamma(L^2(0,T;H), E))$;

(4) *For almost all* $\omega \in \Omega$, *the function* $\Phi_\omega$ *is stochastically integrable with respect to an independent $H$-cylindrical Brownian motion* $\tilde{W}_H$, *and* $\omega \mapsto \int_0^T \Phi(t, \omega) \, d\tilde{W}_H(t)$ *defines an element of* $L^p(\Omega; L^p(\tilde{\Omega}; E))$.

*In this situation the random variables* $\eta$ *in* (1) *and* (2) *are uniquely determined and equal as elements of* $L^p(\Omega; E)$, *the element $X$ in* (3) *is in* $L^p_{\mathbb{F}}(\Omega; \gamma(L^2(0,T;H), E))$, *and we have* $\eta = I^{W_H}(X)$ *in* $L^p(\Omega; E)$. *Moreover,*

$$(3.2) \qquad \mathbb{E} \|X\|^p_{\gamma(L^2(0,T;H), E)} \eqsim_p \mathbb{E} \left\| \int_0^T \Phi(t) \, d\tilde{W}_H(t) \right\|^p_{L^p(\tilde{\Omega}; E)}$$

*and*

$$(3.3) \qquad \beta_{p,E}^{-p} \mathbb{E} \|X\|^p_{\gamma(L^2(0,T;H), E)} \lesssim_p \mathbb{E} \|\eta\|^p \lesssim_p \beta_{p,E}^p \mathbb{E} \|X\|^p_{\gamma(L^2(0,T;H), E)}.$$

A process $\Phi : [0, T] \times \Omega \to \mathcal{L}(H, E)$ satisfying the equivalent conditions of the theorem will be called $L^p$-*stochastically integrable* with respect to $W_H$. The random variable $\eta = I^{W_H}(X)$ is called the *stochastic integral* of $\Phi$ with respect to $W_H$, notation

$$\eta = I^{W_H}(X) =: \int_0^T \Phi(t) \, dW_H(t).$$

REMARK 3.7. Under the assumptions as stated, condition (1) is equivalent to:

(1)′ *There exists a sequence* $(\Phi_n)_{n \geq 1}$ *of elementary adapted processes such that:*

(i) *for all* $h \in H$ *we have* $\lim_{n \to \infty} \Phi_n h = \Phi h$ *in measure on* $[0, T] \times \Omega$;
(ii) *there exists an* $\eta \in L^p(\Omega; E)$ *such that*

$$\eta = \lim_{n \to \infty} \int_0^T \Phi_n(t) \, dW_H(t) \qquad \text{in } L^p(\Omega; E).$$

The proof, as well as further approximation results, will be presented elsewhere.

PROOF OF THEOREM 3.6. (4) $\Leftrightarrow$ (3): This equivalence follows from Lemma 2.7; together with (3.1) this also gives (3.2).

(3) ⇒ (1): By Propositions 2.11 and 2.12, $X \in L^p(\Omega; \gamma(L^2(0,T;H), E))$ represented by $\Phi$ belongs to $L^p_{\mathbb{F}}(\Omega; \gamma(L^2(0,T;H), E))$. Thus we may choose a sequence $(X_n)_{n\geq 1}$ of elementary adapted elements with $\lim_{n\to\infty} X_n = X$ in $L^p(\Omega; \gamma(L^2(0,T;H), E))$. Let $(\Phi_n)_{n\geq 1}$ be a representing sequence of elementary adapted processes. The sequence $(\Phi_n)_{n\geq 1}$ has properties (i) and (ii). Indeed, property (i) follows by noting that $\lim_{n\to\infty} \Phi_n^* x^* = \lim_{n\to\infty} \langle X_n, x^* \rangle = \langle X, x^* \rangle = \Phi^* x^*$ in $L^p(\Omega; L^2(0,T;H))$, and hence in measure on $[0,T] \times \Omega$, for all $x^* \in E^*$. Property (ii), with $\eta = I^{W_H}(X)$, follows from Theorem 3.5, since

$$\lim_{n\to\infty} \int_0^T \Phi_n(t)\,dW_H(t) = \lim_{n\to\infty} I^{W_H}(X_n) = I^{W_H}(X) \qquad \text{in } L^p(\Omega; E).$$

The two-sided estimate (3.3) now follows from Theorem 3.5.

(1) ⇒ (2): This follows from the Burkholder–Davis–Gundy inequalities, which imply that for all $x^* \in E^*$ we have $\lim_{n\to\infty} \Phi_n^* x^* = \Phi^* x^*$ in $L^p(\Omega; L^2(0,T;H))$.

(2) ⇒ (3): This is the technical part of the proof. It simplifies considerably for spaces $E$ having a Schauder basis. To get around such an assumption, we give an approximation argument via quotient maps. We proceed in several steps.

We denote by $B_F$ the closed unit ball of a Banach space $F$.

Since $\Phi$ is $H$-strongly measurable and adapted, without loss of generality we may assume that $E$ is separable. Since $E$ is reflexive, $E^*$ is separable as well and we may fix a dense sequence $(x_n^*)_{n\geq 1}$ in $B_{E^*}$. Define the closed linear subspaces $F_n$ of $E$ by

$$F_n := \bigcap_{i=1}^n \ker(x_i^*).$$

Let $E_n$ be the quotient space $E/F_n$, and let $Q_n : E \to E_n$ be the quotient map. Then $\dim(E_n) < \infty$ and there is a canonical isomorphism $E_n^* \simeq F_n^\perp$, where $F_n^\perp = \{x^* \in E^* : x^* = 0 \text{ on } F_n\}$.

*Step* 1. For every finite-dimensional subspace $G$ of $E$ and every $\varepsilon > 0$ there exists an index $N \geq 1$ such that

(3.4) $$\|x\| \leq (1+\varepsilon)\|Q_N x\| \qquad \forall x \in G.$$

To show this it suffices to consider $x \in B_G$. Since $B_G$ is compact we can find elements $y_1^*, \ldots, y_n^* \in E^*$ with $\|y_i^*\| \leq 1$ such that

$$\|x\| \leq \left(1 + \frac{\varepsilon}{2}\right) \sup_{1\leq i\leq n} |\langle x, y_i^* \rangle| \qquad \forall x \in B_G.$$

Since $(x_i^*)_{i\geq 1}$ is norm dense in $B_{E^*}$, we may approximate the $y_i^*$ to obtain an index $N$ such that

$$\|x\| \leq (1+\varepsilon) \sup_{1\leq j\leq N} |\langle x, x_j^* \rangle| \qquad \forall x \in B_G.$$

It follows that for all $x \in B_G$,
$$\|x\| \leq (1+\varepsilon) \inf_{y \in F_N} \sup_{1 \leq j \leq N} |\langle x - y, x_j^* \rangle| \leq (1+\varepsilon) \inf_{y \in F_N} \|x - y\| = (1+\varepsilon)\|Q_N x\|.$$

This proves (3.4).

*Step* 2. Let the processes $\Phi_n : [0,T] \times \Omega \to \mathcal{L}(H; E_n)$ be given by $\Phi_n(t,\omega)h := Q_n \Phi(t,\omega)h$. Clearly $\Phi_n$ belongs to $L^p(\Omega; L^2(0,T;H))$ scalarly. Moreover, $\Phi_n$ represents an element $X_n \in L^p(\Omega; \gamma(L^2(0,T;H), E_n))$, since for the finite-dimensional spaces $E_n$ we have $\gamma(L^2(0,T;H), E_n) \simeq \mathcal{L}(L^2(0,T;H), E_n)$. Note that almost surely, in $L^2(0,T;H)$ we have

(3.5) $\qquad \langle X_n, x^* \rangle = \Phi_n^* x^* \qquad$ for all $x^* \in E^*$.

This can be proved directly or deduced from Lemma 2.7.

It is easily checked that $I^{W_H} X_n = Q_n \eta$. Hence,
$$\mathbb{E}\|X_n\|_{\gamma(L^2(0,T;H),E_n)}^p \lesssim_p \beta_{p,E_n}^p \mathbb{E}\|I^{W_H} X_n\|_{E_n}^p = \beta_{p,E_n}^p \mathbb{E}\|Q_n \eta\|_{E_n}^p$$
$$\stackrel{(*)}{\leq} \beta_{p,E}^p \mathbb{E}\|Q_n \eta\|_{E_n}^p \leq \beta_{p,E}^p \mathbb{E}\|\eta\|^p.$$

In $(*)$ we used the well known fact that the UMD$(p)$ constant of a quotient space of $E$ can be estimated by the UMD$(p)$ constant of $E$.

For $1 \leq m \leq n$ let $Q_{nm} : E_n \to E_m$ be given by $Q_{nm} Q_n x := Q_m x$. Then $\|Q_{nm}\| \leq 1$ and $X_m = Q_{nm} X_n$. It follows that $\mathbb{E}\|X_m\|_{\gamma(L^2(0,T;H),E_m)} \leq \mathbb{E}\|X_n\|_{\gamma(L^2(0,T;H),E_n)}$. By Fatou's lemma,

(3.6) $\quad \mathbb{E} \sup_{n \geq 1} \|X_n\|_{\gamma(L^2(0,T;H),E_n)}^p = \mathbb{E} \lim_{n \to \infty} \|X_n\|_{\gamma(L^2(0,T;H),E_n)}^p \lesssim_{p,E} \mathbb{E}\|\eta\|^p.$

*Step* 3. Let $N_0$ be a null set such that for all $\omega \in \complement N_0$ we have
$$C(\omega) := \sup_{n \geq 1} \|X_n(\omega)\|_{\gamma(L^2(0,T;H),E_n)} < \infty.$$

Using (3.5), for each $n \geq 1$ we can find a null set $N_n$ of that for all $\omega \in \complement N_n$ and $x^* \in E_n^*$, $\langle X_n(\omega), x^* \rangle = \Phi_n^*(\cdot, \omega)x^*$ in $L^2(0,T;H)$. Let $N := N_0 \cup (\bigcup_{n \geq 1} N_n)$. We claim that for all $\omega \in \complement N$ and all $x^* \in E^*$, $\Phi^*(\cdot, \omega)x^* \in L^2(0,T;H)$.

Fix $\omega \in \complement N$. First let $x^*$ be a linear combination of the elements $x_1^*, \ldots, x_n^*$. Then $x^* \in F_n^\perp$ and hence, for all $t \in [0,T]$, $\Phi^*(t,\omega)x^* = \Phi_n^*(t,\omega)x^*$. It follows that
$$\|\Phi^*(\cdot, \omega)x^*\|_{L^2(0,T;H)} = \|\langle X_n(\omega), x^* \rangle\|_{L^2(0,T;H)}$$
$$\leq \|X_n(\omega)\|_{\gamma(L^2(0,T;H),E_n)} \|x^*\| \leq C(\omega)\|x^*\|.$$

Next let $x^* \in E^*$ be arbitrary; we may assume that $x^* \in B_{E^*}$. Since $(x_k^*)_{k \geq 1}$ is norm dense in $B_{E^*}$ we can find a subsequence $(k_n)_{n \geq 1}$ such that $x^* = \lim_{n \to \infty} x_{k_n}^*$ strongly. It follows that for all $m, n \geq 1$ we have
$$\|\Phi^*(\cdot, \omega)(x_{k_n}^* - x_{k_m}^*)\|_{L^2(0,T;H)} \leq C(\omega)\|x_{k_n}^* - x_{k_m}^*\|.$$

We deduce that $(\Phi^*(\cdot,\omega)x^*_{k_n})_{n\geq 1}$ is a Cauchy sequence in $L^2(0,T;H)$, and after passing to an almost everywhere convergent limit we find that the limit equals $\Phi^*(\cdot,\omega)x^*$. Hence, $\Phi^*(\cdot,\omega)x^* = \lim_{n\to\infty}\Phi^*(\cdot,\omega)x^*_{k_n}$ in $L^2(0,T;H)$. Since $\omega \in \complement N$ was arbitrary, this proves the claim.

*Step* 4. By Step 3, for $\omega \in \complement N$ fixed we may define the integral operator $X(\omega): L^2(0,T;H) \to E$ by

$$X(\omega)f := \int_0^T \Phi(t,\omega)f(t)\,dt.$$

These integrals are well defined as Pettis integrals in $E$ since $E$ is reflexive. We claim that $X(\omega) \in \gamma(L^2(0,T;H),E)$ and

(3.7) $\qquad \|X(\omega)\|_{\gamma(L^2(0,T;H),E)} \leq \sup_{n\geq 1}\|X_n(\omega)\|_{\gamma(L^2(0,T;H),E_n)}.$

To prove this, let the random variables $\rho_n(\omega) \in L^p(\Omega';E)$ be given by

$$\rho_n(\omega) := \sum_{i=1}^n \gamma_i \int_0^T \Phi(t,\omega)f_i(t)\,dt,$$

where $(\gamma_i)_{i\geq 1}$ is a standard Gaussian sequence defined on a probability space $(\Omega',\mathcal{F}',\mathbb{P}')$ and $(f_i)_{i\geq 1}$ is an orthonormal basis for $L^2(0,T;H)$.

Let $\varepsilon > 0$ be arbitrary and fixed. Since $\rho_n(\omega)$ takes its values in a finite-dimensional subspace of $E$, it follows from Step 1 that there is an index $N_n$ such that

$$\mathbb{E}'\|\rho_n(\omega)\|^2 \leq (1+\varepsilon)^2 \mathbb{E}'\|Q_{N_n}\rho_n(\omega)\|^2.$$

Clearly,

$$\mathbb{E}'\|Q_{N_n}\rho_n(\omega)\|^2 = \mathbb{E}'\left\|\sum_{i=1}^n \gamma_i \int_0^T \Phi_{N_n}(t,\omega)f_i(t)\,dt\right\|^2$$
$$\leq \|X_{N_n}(\omega)\|^2_{\gamma(L^2(0,T;H),E_{N_n})},$$

and it follows that

$$\sup_{n\geq 1}\mathbb{E}'\|\rho_n(\omega)\|^2 \leq (1+\varepsilon)^2 \sup_{N\geq 1}\|X_N(\omega)\|^2_{\gamma(L^2(0,T;H),E_N)}.$$

Since $E$ does not contain a copy of $c_0$, a theorem of Hoffmann-Jorgensen and Kwapień [22], Theorem 9.29, assures that $X(\omega) \in \gamma(L^2(0,T;H),E)$ and

$$\|X(\omega)\|^2_{\gamma(L^2(0,T;H),E)} = \sup_{n\geq 1}\mathbb{E}'\|\rho_n(\omega)\|^2 \leq (1+\varepsilon)^2 \sup_{N\geq 1}\|X_N(\omega)\|^2_{\gamma(L^2(0,T;H),E_N)}.$$

Since $\varepsilon > 0$ was arbitrary, the claim follows.

*Step* 5. To finish the proof, we note that $X: \Omega \to \gamma(L^2(0,T;H),E)$ is almost surely equal to a strongly measurable random variable; see Remark 2.8. It follows from (3.6) and (3.7) that $X \in L^p(\Omega;\gamma(L^2(0,T;H),E))$. By definition $X$ is represented by $\Phi$ and hence (3) follows. $\square$

REMARK 3.8. If the filtration $\mathbb{F}$ is assumed to be the augmented Brownian filtration $\mathbb{F}^{W_H}$, the equivalence $(1) \Leftrightarrow (2)$ is true for arbitrary real Banach spaces $E$. This follows from the martingale representation theorem in finite dimensions. We briefly sketch a proof of $(2) \Rightarrow (1)$. For $K = 1, 2, \ldots$ let $\mathcal{F}_T^{(K)}$ be the $\sigma$-algebra generated by the Brownian motions $W_H h_k$, $1 \le k \le K$. Choose a sequence of simple random variables $(\eta_n)_{n \ge 1}$ in $L^p(\Omega, \mathcal{F}_T^{(K)}; E)$ with mean zero and such that $\eta = \lim_{n \to \infty} \eta_n$. This is possible by the martingale convergence theorem and the Pettis measurability theorem. By the martingale representation theorem for finite-dimensional spaces, for all $n \ge 1$ there exists an $L^p$-stochastically integrable process $\Phi_n$ such that $\eta_n = \int_0^T \Phi_n(t)\, dW_H(t)$. The sequence $(\Phi_n)_{n \ge 1}$ satisfies (i) and (ii) of condition (1) of Theorem 3.6. Indeed, (ii) is obvious and (i) follows from the Burkholder–Davis–Gundy inequalities. The processes $\Phi_n$ need not be elementary adapted, but since each $\Phi_n$ takes values in a finite dimensional subspace of $E$ one can approximate the $\Phi_n$ with elementary adapted processes to complete the proof.

For $H = \mathbb{R}$, the implication $(4) \Rightarrow (1)$ in Theorem 3.6 can be interpreted as an $L^p$-version of McConnell's result quoted in the Introduction. Below, in the implication $(4) \Rightarrow (1)$ of Theorem 5.9, we recover McConnell's result.

COROLLARY 3.9 (Series expansion). *Let $E$ be a UMD space and fix $p \in (1, \infty)$. Assume that the $H$-strongly measurable and adapted process $\Phi : [0, T] \times \Omega \to \mathcal{L}(H, E)$ is $L^p$-stochastically integrable with respect to $W_H$. Then for all $h \in H$ the process $\Phi h : [0, T] \times \Omega \to E$ is $L^p$-stochastically integrable with respect to $W_H h$. Moreover, if $(h_n)_{n \ge 1}$ is an orthonormal basis for $H$, then*

$$\int_0^T \Phi(t)\, dW_H(t) = \sum_{n \ge 1} \int_0^T \Phi(t) h_n\, dW_H(t) h_n,$$

*with unconditional convergence in $L^p(\Omega; E)$.*

PROOF. Let $P_N$ be the orthogonal projection in $H$ onto the span of the vectors $h_1, \ldots, h_N$. Let $X \in L^p(\Omega; \gamma(L^2(0, T; H), E))$ be the element represented by $\Phi$. By the right ideal property we have

$$\|X \circ P_N\|_{\gamma(L^2(0,T;H),E)} \le \|X\|_{\gamma(L^2(0,T;H),E)}$$

almost surely. Here we think of $P_N$ as an operator on $\gamma(L^2(0, T; H), E)$ defined by $(P_N S)f := S(P_N f)$ with $(P_N f)(t) := P_N(f(t))$. By an approximation argument one can show that

$$\lim_{N \to \infty} \|X - X \circ P_N\|_{\gamma(L^2(0,T;H),E)} = 0,$$

almost surely. Since $\Phi P_N$ is represented by $X \circ P_N$, the result follows from Theorem 3.6 and the dominated convergence theorem. The convergence of the series is unconditional since any permutation of $(h_n)_{n\geq 1}$ is again an orthonormal basis for $H$. $\square$

A theory of stochastic integration for processes in martingale type 2 spaces has been developed by a number of authors, including Belopolskaya and Daletskiĭ [1], Brzeźniak [4, 5, 6], Dettweiler [11, 12], Neidhardt [30] and Ondreját [31]. Some of these authors state their results for 2-uniformly smooth Banach spaces; the equivalence of martingale type 2 and 2-uniform smoothness up to renorming was shown by Pisier [32]. To make the link with our results, first we recall that a UMD space has martingale (co)type 2 if and only if it has (co)type 2, (cf. [6, 33]), and that every space with martingale type 2 has type 2. By the results of [29, 36], $E$ has type 2 if and only if we have an inclusion $L^2(0,T;\gamma(H,E)) \hookrightarrow \gamma(L^2(0,T;H),E)$, and that $E$ has cotype 2 if and only if we have an inclusion $\gamma(L^2(0,T;H),E) \hookrightarrow L^2(0,T;\gamma(H,E))$; in both cases the inclusion is given via representation. Thus from Theorem 3.6 we obtain the following result.

COROLLARY 3.10. *Let $E$ be a UMD space and let $p \in (1,\infty)$.*

(1) *If $E$ has type $2$, then every $H$-strongly measurable and adapted process $\Phi$ which belongs to $L^p(\Omega; L^2(0,T;\gamma(H,E)))$ is $L^p$-stochastically integrable with respect to $W_H$ and we have*

$$\mathbb{E}\left\|\int_0^T \Phi(t)\,dW_H(t)\right\|^p \lesssim_{p,E} \mathbb{E}\|\Phi\|_{L^2(0,T;\gamma(H,E))}^p.$$

(2) *If $E$ has cotype $2$, then every $H$-strongly measurable process $\Phi$ which is $L^p$-stochastically integrable with respect to $W_H$ belong to $L^p(\Omega; L^2(0,T;\gamma(H,E)))$ and we have*

$$\mathbb{E}\|\Phi\|_{L^2(0,T;\gamma(H,E))}^p \lesssim_{p,E} \mathbb{E}\left\|\int_0^T \Phi(t)\,dW_H(t)\right\|^p.$$

We conclude this section with a result giving a necessary and sufficient square function criterion for $L^p$-stochastic integrability of $\mathcal{L}(H,E)$-valued processes, where $E$ is assumed to be a UMD Banach function space. In view of Theorem 3.6 it suffices to give such a criterion for $\mathcal{L}(H,E)$-valued functions, and therefore a straightforward extension of [28], Corollary 2.10 (where only the case $H = \mathbb{R}$ was considered) gives the following result.

COROLLARY 3.11. *Let $E$ be UMD Banach function space over a $\sigma$-finite measure space $(S,\Sigma,\mu)$ and let $p \in (1,\infty)$. Let $\Phi:[0,T] \times \Omega \to \mathcal{L}(H,E)$ be $H$-strongly measurable and adapted and assume that there exists a strongly* measurable function $\phi\colon [0,T] \times \Omega \times S \to H$ such that for all $h \in H$ and $t \in [0,T]$,

$$(\Phi(t)h)(\cdot) = [\phi(t,\cdot), h]_H \qquad \text{in } E.$$

Then $\Phi$ is $L^p$-stochastically integrable with respect to $W_H$ if and only if

$$\mathbb{E}\left\|\left(\int_0^T \|\phi(t,\cdot)\|_H^2\, dt\right)^{1/2}\right\|_E^p < \infty.$$

In this case we have

$$\mathbb{E}\left\|\int_0^T \Phi(t)\, dW_H(t)\right\|^p \eqsim_{p,E} \mathbb{E}\left\|\left(\int_0^T \|\phi(t,\cdot)\|_H^2\, dt\right)^{1/2}\right\|_E^p.$$

**4. The integral process.** It is immediate from Theorem 3.6 that if $\Phi\colon [0,T] \times \Omega \to \mathcal{L}(H,E)$ is $L^p$-stochastically integrable with respect to $W_H$, then for all $t \in [0,T]$ the restricted process $\Phi\colon [0,t] \times \Omega \to \mathcal{L}(H,E)$ is $L^p$-stochastically integrable with respect to $W_H$. Thus it is meaningful to ask for the properties of the integral process

$$t \mapsto \int_0^t \Phi(s)\, dW_H(s), \qquad t \in [0,T].$$

This will be the topic of the present section.

It will be convenient to introduce a continuous process

$$\xi_X \colon [0,T] \times \Omega \to \gamma(L^2(0,T;H), E)$$

associated with a strongly measurable random variable $X\colon \Omega \to \gamma(L^2(0,T; H), E)$. For $t \in [0,T]$ we define the $\gamma(L^2(0,T;H), E)$-valued random variable $\xi_X(t)\colon \Omega \to \gamma(L^2(0,T;H), E)$ by

$$\xi_X(t,\omega)f := (X(\omega))(\mathbf{1}_{[0,t]}f), \qquad f \in L^2(0,T;H).$$

Note that $\xi_X(T) = X$. The strong measurability of $\xi_X(t)$ as a $\gamma(L^2(0,T; H), E)$-valued random variable follows from Lemma 2.5.

PROPOSITION 4.1. *The process $\xi_X$ defined above is strongly measurable and has continuous trajectories. Moreover:*

(1) *If $X$ is strongly adapted to $\mathbb{F}$, then $\xi_X$ is adapted to $\mathbb{F}$ and for all $t \in [0,T]$, $\xi_X(t)$ is strongly adapted to $\mathbb{F}$;*

(2) *If $X \in L^p_{\mathbb{F}}(\Omega; \gamma(L^2(0,T;H), E))$, then $\xi_X(t) \in L^p_{\mathbb{F}}(\Omega; \gamma(L^2(0,T;H), E))$ for all $t \in [0,T]$, and the mapping $t \mapsto \xi_X(t)$ is continuous from $[0,T]$ to $L^p_{\mathbb{F}}(\Omega; \gamma(L^2(0,T;H), E))$.*

PROOF. By Proposition 2.4, $t \mapsto \xi_X(t,\omega)$ is continuous for all $\omega \in \Omega$. Since for all $t \in [0,T]$, $\xi_X(t)$ is strongly measurable we obtain that $\xi_X$ is strongly measurable.

(1) This follows from Lemma 2.5 and Proposition 2.10.
(2) For $\omega \in \Omega$ fixed, the right ideal property implies that

$$\|\xi_X(t)(\omega)\|_{\gamma(L^2(0,T;H),E)} \leq \|X(\omega)\|_{\gamma(L^2(0,T;H),E)}.$$

Hence if $X \in L^p_{\mathbb{F}}(\Omega; \gamma(L^2(0,T;H),E))$, then for all $t \in [0,T]$, $\xi_X(t) \in L^p_{\mathbb{F}}(\Omega; \gamma(L^2(0,T;H),E))$ by Proposition 2.12. The continuity of $t \mapsto \xi_X(t)$ follows from Proposition 4.1 and dominated convergence. $\square$

REMARK 4.2. Since $(t,\omega) \mapsto \|\xi_X(t,\omega)\|^2_{\gamma(L^2(0,T;H),E)}$ is nonnegative and nondecreasing, we may think of this process as an analogue of the quadratic variation process.

Now let $E$ be a UMD space and fix $p \in (1,\infty)$. For $X \in L^p_{\mathbb{F}}(\Omega; \gamma(L^2(0,T;H), E))$, with some abuse of notation the $E$-valued process

$$I^{W_H}(\xi_X) : t \mapsto I^{W_H}(\xi_X(t)), \qquad t \in [0,T],$$

will be called the *integral process* associated with $X$. In the special case where $X$ is represented by an $L^p$-stochastically integrable process $\Phi$, for all $t \in [0,T]$ we have

$$I^{W_H}(\xi_X(t)) = \int_0^t \Phi(s)\, dW_H(s) \qquad \text{in } L^p(\Omega;E).$$

PROPOSITION 4.3. *Let $E$ be a UMD space and fix $p \in (1,\infty)$. For all $X \in L^p_{\mathbb{F}}(\Omega; \gamma(L^2(0,T;H),E))$ the integral process $I^{W_H}(\xi_X)$ is an $E$-valued $L^p$-martingale which is continuous in $p$th moment. It has a continuous adapted version which satisfies the maximal inequality*

(4.1) $\quad \mathbb{E} \sup_{t \in [0,T]} \|I^{W_H}(\xi_X(t))\|^p \leq q^p \mathbb{E} \|I^{W_H}(X)\|^p \qquad \left(\frac{1}{p} + \frac{1}{q} = 1\right).$

PROOF. For all $x^* \in E^*$, the real-valued process $I^{W_H}(\xi^*_X x^*)$ is a martingale; see [18], Corollary 17.8. The martingale property easily follows from this; see [28], Corollary 2.8. The continuity in $p$th moment follows directly from the continuity of the Itô map and the continuity in $p$th moment of $\xi_X$, which was proved in Proposition 4.1.

Next we prove the existence of a continuous adapted version. Choose a sequence $(X_n)_{n \geq 1}$ of elementary adapted elements such that $\lim_{n \to \infty} X_n = X$ in $L^p(\Omega; \gamma(L^2(0,T;H),E))$. It follows from Theorem 3.5 that $\lim_{n \to \infty} I^{W_H}(X_n) = I^{W_H}(X)$ in $L^p(\Omega; E)$. Clearly, for each $n \geq 1$ there exists a continuous version $\eta_n$ of $I^{W_H}(\xi^{X_n})$, and by the Pettis measurability theorem we have $\eta_n \in L^p(\Omega; C([0,T]; E))$. By Doob's maximal inequality, the sequence $(\eta_n)_{n \geq 1}$ is a Cauchy sequence in $L^p(\Omega; C([0,T]; E))$. Its limit defines a continuous version of $I^{W_H}(\xi_X)$, which is clearly adapted.

The final inequality (4.1) follows from Doob's maximal inequality. □

Combining these results we have proved:

THEOREM 4.4 (Burkholder–Davis–Gundy inequalities). *Let $E$ be a UMD space and fix $p \in (1, \infty)$. If the $H$-strongly measurable and adapted process $\Phi:[0,T] \times \Omega \to \mathcal{L}(H,E)$ is $L^p$-stochastically integrable, then*

$$\mathbb{E} \sup_{t \in [0,T]} \left\| \int_0^t \Phi(s)\, dW_H(s) \right\|^p \eqsim_{p,E} \mathbb{E}\|X\|^p_{\gamma(L^2(0,T;H),E)},$$

*where $X \in L^p(\Omega; \gamma(L^2(0,T;H),E))$ is the element represented by $\Phi$.*

The estimates in Corollary 3.10, when combined with Doob's maximal inequality, may be considered as one-sided Burkholder–Davis–Gundy inequalities for the $L^p(\Omega; L^2(0,T; \gamma(H,E)))$-norm. In particular we recover, for UMD martingale type 2 spaces, the one-sided Burkholder–Davis–Gundy inequalities for martingale type 2 spaces of Brzeźniak [6] and Dettweiler [12].

We address next the question whether the integral process associated with an $L^p$-stochastically integrable process $\Phi$ is $L^p$-stochastically integrable with respect to a real-valued Brownian motion $W$. When $E$ is a real Hilbert space and $p \in (1, \infty)$, the answer is clearly affirmative and by the Burkholder–Davis–Gundy inequalities we have

$$\left( \mathbb{E} \left\| \int_0^T \int_0^t \Phi(s)\, dW_H(s)\, dW(t) \right\|^p \right)^{1/p}$$
$$\eqsim_p \left\| \int_0^\cdot \Phi(s)\, dW_H(s) \right\|_{L^p(\Omega; L^2(0,T;E))}$$
$$\leq \sqrt{T} \left( \mathbb{E} \sup_{t \in [0,T]} \left\| \int_0^t \Phi(s)\, dW_H(s) \right\|^p \right)^{1/p} \eqsim_p \sqrt{T} \|\Phi\|_{L^p(\Omega; L^2(0,T;E))}.$$

More generally, every $\mathcal{L}_2(H,E)$-valued $L^p$-martingale, where $E$ is a Hilbert space, is $L^p$-stochastically integrable, and an estimate can be given using Doob's inequality. In the following we shall generalize these observations to $\gamma(H,E)$-valued $L^p$-martingales, where $E$ is a UMD space. We will say that a process $M:[0,T] \times \Omega \to \gamma(H,E)$ is an $L^p$-*martingale* if $M(t) \in L^p(\Omega; \gamma(H,E))$ for all $t \in [0,T]$ and $\mathbb{E}(M(t)|\mathcal{F}_s) = M(s)$ in $L^p(\Omega; \gamma(H,E))$ for all $0 \leq s \leq t \leq T$. In the proof of the following result we will need the well known fact that every $L^p$-martingale $M:[0,T] \times \Omega \to H$ admits a modification with cadlag trajectories. This may be proved as [21], Proposition 2.

Our next result uses the vector-valued Stein inequality, which asserts that in a UMD space $E$ certain families of conditional expectation operators are $R$-bounded. Recall that a collection $\mathcal{T} \subseteq \mathcal{L}(B_1, B_2)$, where $B_1$ and $B_2$ are Banach spaces, is said to be *R-bounded* if there exists a constant $M \geq 0$ such that

$$\left( \mathbb{E} \left\| \sum_{n=1}^{N} r_n T_n x_n \right\|_{B_2}^2 \right)^{1/2} \leq M \left( \mathbb{E} \left\| \sum_{n=1}^{N} r_n x_n \right\|_{B_1}^2 \right)^{1/2},$$

for all $N \geq 1$ and all sequences $(T_n)_{n=1}^N$ in $\mathcal{T}$ and $(x_n)_{n=1}^N$ in $B_1$. The least constant $M$ for which this estimate holds is called the *R-bound* of $\mathcal{T}$, notation $R(\mathcal{T})$. By the Kahane–Khinchine inequalities, the role of the exponent 2 may be replaced by any exponent $1 \leq p < \infty$. Replacing the role of the Rademacher sequence by a Gaussian sequence we obtain the related notion of $\gamma$-boundedness. By an easy randomization argument, every $R$-bounded family is $\gamma$-bounded and we have $\gamma(\mathcal{T}) \leq R(\mathcal{T})$, where $\gamma(\mathcal{T})$ is the *$\gamma$-bound* of $\mathcal{T}$.

THEOREM 4.5. *Let $E$ be a UMD space and fix $p \in (1, \infty)$. Let $M: [0, T] \times \Omega \to \gamma(H, E)$ be an $L^p$-martingale with respect to the filtration $\mathbb{F}$ and assume that $M(0) = 0$. If $W_H$ is an $H$-cylindrical Brownian motion adapted to $\mathbb{F}$, then $M$ is $L^p$-stochastically integrable with respect to $W_H$ and we have*

$$\left( \mathbb{E} \left\| \int_0^T M(t) \, dW_H(t) \right\|^p \right)^{1/p} \lesssim_{p,E} \sqrt{T} (\mathbb{E} \| M(T) \|_{\gamma(H,E)}^p)^{1/p}.$$

PROOF. The proof is based upon a multiplier result for spaces of $\gamma$-radonifying operators, due to Kalton and the third named author [19]. Translated into the present setting, this result can be formulated as follows. Let $B_1$ and $B_2$ be UMD spaces, let $p \in (1, \infty)$, and let $N: [0, T] \times \Omega \to \mathcal{L}(B_1, B_2)$ be a strongly adapted process such that the set $\{N(t): t \in [0, T]\}$ is $\gamma$-bounded. Then, if $\Phi: [0, T] \times \Omega \to \mathcal{L}(H, B_1)$ is an $H$-strongly measurable process which is $L^p$-stochastically integrable with respect to $W_H$, the process $N\Phi: [0, T] \times \Omega \to \mathcal{L}(H, B_2)$ defined by $(N\Phi)(t)h := N(t)(\Phi(t)h)$ is $L^p$-stochastically integrable with respect to $W_H$ as well and satisfies

$$\mathbb{E} \left\| \int_0^T N(t) \Phi(t) \, dW_H(t) \right\|^p \lesssim_{p, B_1, B_2} K^p \mathbb{E} \left\| \int_0^T \Phi(t) \, dW_H(t) \right\|^p.$$

To start the proof of the theorem, we first show that $M$ is $H$-strongly measurable and adapted. Let $h \in H$ be fixed. Clearly, $Mh$ is an $E$-valued $L^p$-martingale. By martingale convergence, $Mh$ is left continuous in mean. Therefore by a general result from the theory of stochastic processes, $Mh$ is strongly measurable and adapted.

Next we check that $M$ belongs to $L^p(\Omega; L^2(0, T; H))$ scalarly. Let $x^* \in E^*$ be fixed. By the above discussion $M^* x^*$ has a modification with cadlag trajectories. Hence we may apply Doob's maximal inequality to obtain

$$\mathbb{E}\|M^*x^*\|^p_{L^2(0,T;H)} \leq T^{p/2}\mathbb{E}\sup_{t\in[0,T]}\|M^*(t)x^*\|^p_H \lesssim_p T^{p/2}\mathbb{E}\|M^*(T)x^*\|^p_H.$$

Let $B = L^p_0(\Omega, \mathcal{F}_T; E)$ be the closed subspace in $L^p(\Omega; E)$ of all $\mathcal{F}_T$-measurable random variables with zero mean, and define the bounded and strongly left continuous function $N:[0,T] \to \mathcal{L}(B)$ by

$$N(t)\xi := \mathbb{E}(\xi|\mathcal{F}_t), \qquad \xi \in B, t \in [0,T].$$

Since $E$ is a UMD space, by a result of Bourgain [3] the set $\{N(t) : t \in [0,T]\}$ is $R$-bounded, and therefore $\gamma$-bounded, with $\gamma$-bound depending only on $p$ and $E$. A detailed proof of this fact may be found in [10], Proposition 3.8.

By the Fubini isomorphism we may identify the random variables $M(t) \in L^p(\Omega; \gamma(H, E))$ with operators $\tilde{M}(t) \in \gamma(H, L^p(\Omega; E))$. Recall that for all $t \in [0,T]$, for all $h \in H$, for almost all $\omega \in \Omega$, $(\tilde{M}(t)h)(\omega) = M(t,\omega)h$. Define a constant function $G:[0,T] \to \mathcal{L}(H, B)$ by

$$G(t) := \tilde{M}(T), \qquad t \in [0,T].$$

Clearly $G$ represents the element $R_G \in \gamma(L^2(0,T;H), B)$ given by

$$R_G f = \int_0^T \tilde{M}(T)f(t)\,dt, \qquad f \in L^2(0,T;H),$$

and $\|R_G\|_{\gamma(L^2(0,T;H),B)} = \sqrt{T}\mathbb{E}\|M(T)\|_{\gamma(H,E)}$. Since for all $t \in [0,T]$, $\tilde{M}(t) = N(t)\tilde{M}(T)$ in $B$, we may apply the above multiplier result to conclude that $\tilde{M}$ represents an element $R \in \gamma(L^2(0,T;H), B)$ with

$$\|R\|_{\gamma(L^2(0,T;H),B)} \lesssim_{p,E} \|R_G\|_{\gamma(L^2(0,T;H),B)}.$$

Using the $\gamma$-Fubini isomorphism we define $X = F_\gamma^{-1}(R)$. Recall that for all $f \in L^2(0,T;H)$, for almost all $\omega \in \Omega$, $(Rf)(\omega) = X(\omega)f$.

We claim that $X$ is represented by $M$. Once we know this, it follows with Theorem 3.6 that

$$\left(\mathbb{E}\left\|\int_0^T M(t)\,dW_H(t)\right\|^p\right)^{1/p} \eqsim_{p,E} \left(\mathbb{E}\|X\|^p_{\gamma(L^2(0,T;H),E)}\right)^{1/p}$$

$$\eqsim_p \|R\|_{\gamma(L^2(0,T;H),B)}$$

$$\lesssim_{p,E} \sqrt{T}\left(\mathbb{E}\|M(T)\|^p_{\gamma(H,E)}\right)^{1/p}.$$

Let $f \in L^2(0,T;H)$, $x^* \in E^*$ be arbitrary. We have to show that $[M^*x^*, f]_{L^2(0,T;H)} = \langle Xf, x^* \rangle$ almost surely. It suffices to check that $\mathbb{E}(\mathbf{1}_A[M^*x^*,$

$f]_{L^2(0,T;H)}) = \mathbb{E}(\mathbf{1}_A \langle Xf, x^* \rangle)$ for all $A \in \mathcal{F}_T$. By the Fubini theorem we have

$$\mathbb{E}(\mathbf{1}_A[M^*x^*, f]_{L^2(0,T;H)}) = \int_\Omega \int_0^T \langle M(t,\omega)f(t), x^* \rangle \mathbf{1}_A(\omega)\, dt\, dP(\omega)$$
$$= \int_0^T \int_\Omega \langle M(t,\omega)f(t), x^* \rangle \mathbf{1}_A(\omega)\, dP(\omega)\, dt$$
$$= \int_0^T \langle \tilde{M}(t)f(t), \mathbf{1}_A \otimes x^* \rangle\, dt = \langle Rf, \mathbf{1}_A \otimes x^* \rangle$$
$$= \mathbb{E}(\langle Xf, x^* \rangle \mathbf{1}_A).$$

This proves the claim.  □

In view of Proposition 4.3, this theorem can be applied to the integral process $I^{W_H}(\xi_X)$ associated with elements $X \in L^p_\mathbb{F}(\Omega; \gamma(L^2(0,T;H),E))$. In the special case where $X$ is represented by a process we obtain:

COROLLARY 4.6.  *Let $E$ be a UMD space and fix $p \in (1,\infty)$. Let $W_H$ and $W$ be an $H$-cylindrical Brownian motion and a Brownian motion, respectively, both adapted to the filtration $\mathbb{F}$. If the $H$-strongly measurable and adapted process $\Phi:[0,T] \times \Omega \to \mathcal{L}(H,E)$ is $L^p$-stochastically integrable with respect $W_H$, then the integral process $(\int_0^t \Phi(s)\, dW_H(s))_{t \in [0,T]}$ is $L^p$-stochastically integrable with respect to $W$ and we have*

$$\left(\mathbb{E}\left\|\int_0^T \int_0^t \Phi(s)\, dW_H(s)\, dW(t)\right\|^p\right)^{1/p} \lesssim_{p,E} \sqrt{T} \left(\mathbb{E}\left\|\int_0^T \Phi(t)\, dW_H(t)\right\|^p\right)^{1/p}.$$

We conclude this section with a representation theorem for $E$-valued Brownian $L^p$-martingales, that is, $E$-valued $L^p$-martingales adapted to the augmented filtration $\mathbb{F}^{W_H}$ generated by an $H$-cylindrical Brownian motion $W_H$. It is a direct consequence of the second part of Theorem 3.5 and Proposition 4.3:

THEOREM 4.7 (Representation of Brownian $L^p$-martingales in UMD spaces).  *Let $E$ be a UMD space and fix $p \in (1,\infty)$. Then every $L^p$-martingale $M:[0,T] \times \Omega \to E$ adapted to the augmented filtration $\mathbb{F}^{W_H}$ has a continuous version, and there exists a unique $X \in L^p_\mathbb{F}(\Omega; \gamma(L^2(0,T;H),E))$ such that for all $t \in [0,T]$ we have*

$$M(t) = M(0) + I^{W_H}(\xi_X(t))  \quad \text{in } L^p(\Omega;E).$$

**5. Localization.**  We begin with a lemma which is a slight generalization of a stopping time argument in [24], Lemma 3.3. For the convenience of the reader we include the details.

LEMMA 5.1. *Let $p \in [1, \infty)$. Let $E$ and $F$ be Banach spaces and let $(\phi_t)_{t \in [0,T]}$ and $(\psi_t)_{t \in [0,T]}$ be continuous adapted processes with values in $E$ and $F$, respectively. Assume furthermore that $\psi_0 = 0$. If there exists a constant $C \geq 0$ such that for all stopping times $\tau$ with values in $[0, T]$ we have*

$$(5.1) \qquad \mathbb{E}\|\phi_\tau\|_E^p \leq C \mathbb{E}\|\psi_\tau\|_F^p$$

*whenever these norms are finite, then for all $\delta > 0$ and $\varepsilon > 0$ we have*

$$(5.2) \qquad \mathbb{P}\bigg(\sup_{t \in [0,T]} \|\phi_t\|_E > \varepsilon\bigg) \leq \frac{C \delta^p}{\varepsilon^p} + \mathbb{P}\bigg(\sup_{t \in [0,T]} \|\psi_t\|_F \geq \delta\bigg).$$

PROOF. Let $\delta, \varepsilon > 0$ be arbitrary. Define stopping times $\mu$ and $\nu$ by

$$\mu(\omega) := \inf\{t \in [0,T] : \|\phi_t(\omega)\|_E \geq \varepsilon\},$$
$$\nu(\omega) := \inf\{t \in [0,T] : \|\psi_t(\omega)\|_F \geq \delta\},$$

where we take $\mu(\omega) := T$ and $\nu(\omega) := T$ if the infimum is taken over the empty set, and put $\tau := \mu \wedge \nu$. Then $\tau$ is a stopping time and $\mathbb{E}\|\phi_\tau\|_E^p \leq \varepsilon^p$, $\mathbb{E}\|\psi_\tau\|_F^p \leq \delta^p$. By Chebyshev's inequality, (5.1), and pathwise continuity we have

$$\mathbb{P}\bigg(\sup_{t \in [0,T]} \|\phi_t\|_E > \varepsilon, \sup_{t \in [0,T]} \|\psi_t\|_F < \delta\bigg) \leq \mathbb{P}(\|\phi_\tau\|_E \geq \varepsilon) \leq \frac{1}{\varepsilon^p} \mathbb{E}\|\phi_\tau\|_E^p$$
$$\leq \frac{C}{\varepsilon^p} \mathbb{E}\|\psi_\tau\|_F^p \leq \frac{C \delta^p}{\varepsilon^p},$$

where the last inequality uses the fact that $\psi_0 = 0$. This implies

$$\mathbb{P}\bigg(\sup_{t \in [0,T]} \|\phi_t\|_E > \varepsilon\bigg) \leq \frac{C \delta^p}{\varepsilon^p} + \mathbb{P}\bigg(\sup_{t \in [0,T]} \|\phi_t\|_E > \varepsilon, \sup_{t \in [0,T]} \|\psi_t\|_F \geq \delta\bigg).$$

Clearly (5.2) follows from this. □

For a Banach space $B$, let $L^0(\Omega; B)$ be the vector space of all equivalence classes of strongly measurable functions on $\Omega$ with values in the Banach space $B$ which are identical almost surely. Endowed with the translation invariant metric

$$\|\xi\|_{L^0(\Omega; B)} = \mathbb{E}(\|\xi\| \wedge 1),$$

$L^0(\Omega; B)$ is a complete metric space, and convergence with respect to this metric coincides with convergence in probability.

We return to the standing assumptions that $H$ is a separable real Hilbert space, $W_H$ is an $H$-cylindrical Brownian motion adapted to a filtration $\mathbb{F}$ satisfying the usual conditions, and $E$ is a real Banach space. We denote by $L^0_\mathbb{F}(\Omega; \gamma(L^2(0,T;H),E))$ the subspace of all adapted elements of $L^0(\Omega; \gamma(L^2(0,T;H),E))$, that is, the closure of subspace of all elementary adapted elements in $L^0(\Omega; \gamma(L^2(0,T;H),E))$. Notice that $X \in L^0_{\mathbb{F}}(\Omega; \gamma(L^2(0,T;H),E))$ if and only if $X$ is strongly adapted to $\mathbb{F}$.

For a stopping time $\tau$ with values in $[0,T]$ and an element $X \in L^0_{\mathbb{F}}(\Omega; \gamma(L^2(0,T;H),E))$ we define the $\gamma(L^2(0,T;H),E)$-valued random variable $\xi_X(\tau) : \Omega \to \gamma(L^2(0,T;H),E)$ by

$$(\xi_X(\tau))(\omega)f := \xi_X(\tau(\omega),\omega)f = X(\omega)(\mathbf{1}_{[0,\tau(\omega)]}f), \qquad f \in L^2(0,T;H).$$

The random variable $\xi_X(\tau)$ is well-defined since $\xi_X$ has continuous paths and is adapted by Proposition 4.1.

LEMMA 5.2. *The random variable $\xi_X(\tau)$ is strongly adapted to $\mathbb{F}$. If $p \in [1,\infty)$ and $X \in L^p_{\mathbb{F}}(\Omega; \gamma(L^2(0,T;H),E))$, then $\xi_X(\tau) \in L^p_{\mathbb{F}}(\Omega; \gamma(L^2(0,T;H),E))$.*

PROOF. It is clear that for all $t \in [0,T]$, $f \in L^2(0,T;H)$, and $x^* \in E^*$, the random variable $\langle X(\mathbf{1}_{[0,t]}f), x^* \rangle$ is $\mathcal{F}_t$-measurable. Hence the first assertion follows by combining by the Pettis measurability theorem and Proposition 2.10.

By the right ideal property,

$$\|\xi_X(\tau)(\omega)\|_{\gamma(L^2(0,T;H),E)} \leq \|X(\omega)\|_{\gamma(L^2(0,T;H),E)}.$$

Hence if $X \in L^p_{\mathbb{F}}(\Omega; \gamma(L^2(0,T;H),E))$ for some $p \in [1,\infty)$, then $\xi_X(\tau) \in L^p(\Omega; \gamma(L^2(0,T;H),E))$. The second assertion now follows from Proposition 2.12. □

PROPOSITION 5.3. *Let $E$ be a UMD space and let $p \in (1,\infty)$. If $X \in L^p_{\mathbb{F}}(\Omega; \gamma(L^2(0,T;H),E))$ and $\tau$ is a stopping time with values in $[0,T]$, then*

(5.3) $\qquad I^{W_H}(\xi_X(\tau)) = (I^{W_H}(\xi_X))_\tau \qquad$ *almost surely.*

PROOF. For elementary adapted $X$, (5.3) is obvious. For general $X \in L^p_{\mathbb{F}}(\Omega; \gamma(L^2(0,T;H),E))$ the result is obtain from the following approximation argument. Choose a sequence of elementary adapted elements such that $\lim_{n\to\infty} X_n = X$ in $L^p_{\mathbb{F}}(\Omega; \gamma(L^2(0,T;H),E))$. Hence, $\xi_X(\tau) = \lim_{n\to\infty} \xi_{X_n}(\tau)$ in $L^p_{\mathbb{F}}(\Omega; \gamma(L^2(0,T;H),E))$ and it follows from Theorem 3.5 that $I^{W_H}(\xi_X(\tau)) = \lim_{n\to\infty} I^{W_H}(\xi_{X_n}(\tau))$ in $L^p(\Omega; E)$. On the other hand, Proposition 4.3 shows that $I^{W_H}(\xi_X) = \lim_{n\to\infty} I^{W_H}(\xi_{X_n})$ in $L^p(\Omega; C([0,T];E))$. In particular, $(I^{W_H}(\xi_X))_\tau = \lim_{n\to\infty} (I^{W_H}(\xi_{X_n}))_\tau$ in $L^p(\Omega; E)$. The general case of (5.3) now follows from the fact that (5.3) holds for all $X_n$. □

By combining the previous two results we obtain the following result, which should be compared with [24], Lemma 3.3. Our approach is somewhat simpler, as it allows the use of $\mathbb{F}$-stopping times rather than the $\mathbb{F} \otimes \tilde{\mathbb{F}}$-stopping times used in [24].

LEMMA 5.4. *Let $E$ be a UMD space and let $p \in (1, \infty)$. If $X \in L^p_\mathbb{F}(\Omega; \gamma(L^2(0,T;H), E))$, then for all $\delta > 0$ and $\varepsilon > 0$ we have*

$$\text{(5.4)} \quad \mathbb{P}\left(\sup_{t \in [0,T]} \|(I^{W_H}(\xi_X))_t\| > \varepsilon\right) \leq \frac{C_{p,E}\delta^p}{\varepsilon^p} + \mathbb{P}(\|X\|_{\gamma(L^2(0,T;H),E)} \geq \delta)$$

*and*

$$\text{(5.5)} \quad \mathbb{P}(\|X\|_{\gamma(L^2(0,T;H),E)} > \varepsilon) \leq \frac{C_{p,E}\delta^p}{\varepsilon^p} + \mathbb{P}\left(\sup_{t \in [0,T]} \|(I^{W_H}(\xi_X))_t\| \geq \delta\right),$$

*where $C_{p,E}$ is a constant which depends only on $p$ and $E$.*

PROOF. For all $\omega \in \Omega$ and $t \in [0, T]$,

$$\|(\xi_X(t))(\omega)\|_{\gamma(L^2(0,T;H),E)} \leq \|X(\omega)\|_{\gamma(L^2(0,T;H),E)}$$

with equality for $t = T$, and therefore,

$$\|X(\omega)\|_{\gamma(L^2(0,T;H),E)} = \sup_{t \in [0,T]} \|(\xi_X(t))(\omega)\|_{\gamma(L^2(0,T;H),E)}.$$

Hence by Lemma 5.1 it suffices to prove that for every stopping time $\tau$ with values in $[0, T]$ we have

$$\mathbb{E}\|(I^{W_H}(\xi_X))_\tau\|^p \eqsim_{p,E} \mathbb{E}\|\xi_X(\tau)\|^p_{\gamma(L^2(0,T;H),E)}$$

provided both norms are finite. But this follows from Proposition 5.3 and Theorem 3.5. □

We call an $E$-valued process $M := (M_t)_{t \in [0,T]}$ a *local martingale* if it is adapted and there exists a sequence of stopping times $(\tau_n)_{n \geq 1}$ with values in $[0, T]$ with the property that for all $\omega \in \Omega$ there exists an index $N(\omega)$ such that $\tau_n(\omega) = T$ for all $n \geq N(\omega)$ and such that the process $M^{\tau_n} = (M_t^{\tau_n})_{t \in [0,T]}$ defined by

$$M_t^{\tau_n} := M_{t \wedge \tau_n} - M_0$$

is a martingale. In this case, $(\tau_n)_{n \geq 1}$ is called a *localizing sequence* for $M$.

If, for some $p \in [1, \infty]$, each $M^{\tau_n}$ is an $L^p$-martingale, we call $M$ a *local $L^p$-martingale*. In the case of $p = \infty$ we say that $M$ is a *local bounded martingale*. It is easy to see that every continuous local martingale is a continuous local bounded martingale (cf. [9], Proposition 1.9); a localizing sequence $(\tau_n)_{n \geq 1}$ is given by

$$\tau_n = \inf\{t \in [0,T] : \|M_t\| \geq n\}.$$

Here we take $\tau_n = T$ if the infimum is taken over the empty set. We will use this convention for all stopping times in the rest of paper.

We denote by $\mathcal{M}_0^{c,\mathrm{loc}}(\Omega; E)$ the space of continuous local martingales starting at 0, identifying martingales that are indistinguishable. Each $M \in \mathcal{M}_0^{c,\mathrm{loc}}(\Omega; E)$ defines a random variable with values in $C([0,T]; E)$. Thus we may identify $\mathcal{M}_0^{c,\mathrm{loc}}(\Omega; E)$ with a linear subspace of $L^0(\Omega; C([0,T]; E))$. If we want to stress the role of the underlying filtration $\mathbb{F}$ we write $\mathcal{M}_0^{c,\mathrm{loc}}(\Omega; E) = \mathcal{M}_0^{c,\mathrm{loc}}(\Omega, \mathbb{F}; E)$.

Now let $E$ be a UMD space and $p \in (1, \infty)$. For $X \in L_{\mathbb{F}}^p(\Omega; \gamma(L^2(0,T;H), E))$ we recall that from Proposition 4.3 that $I^{W_H}(\xi_X)$ is a continuous martingale starting at 0. With this in mind we have the following localized version of Theorem 3.5.

THEOREM 5.5 (Itô homeomorphism). *Let $E$ be a real UMD space. The mapping $X \mapsto I^{W_H}(\xi_X)$ has a unique extension to a homeomorphism from $L_{\mathbb{F}}^0(\Omega; \gamma(L^2(0,T;H), E))$ onto a closed subspace of $\mathcal{M}_0^{c,\mathrm{loc}}(\Omega, \mathbb{F}; E)$. Moreover, the estimates (5.4) and (5.5) extend to arbitrary elements $X \in L_{\mathbb{F}}^0(\Omega; \gamma(L^2(0,T;H), E))$. For the augmented Brownian filtration $\mathbb{F}^{W_H}$ we have an homeomorphism*

$$I^{W_H} : L_{\mathbb{F}^{W_H}}^0(\Omega; \gamma(L^2(0,T;H), E)) \eqsim \mathcal{M}_0^{c,\mathrm{loc}}(\Omega, \mathbb{F}^{W_H}; E).$$

PROOF. Fix $X \in L_{\mathbb{F}}^0(\Omega; \gamma(L^2(0,T;H), E))$ and define a sequence of stopping times $(\tau_n)_{n \geq 1}$ by

$$\tau_n := \inf\{t \in [0,T] : \|\xi_X(t)\|_{\gamma(L^2(0,T;H),E)} \geq n\}.$$

Then $\xi_X(\tau_n) \in L_{\mathbb{F}}^p(\Omega; \gamma(L^2(0,T;H), E))$ for every $p \in (1, \infty)$.

By Proposition 4.3 we can define a sequence of $L^p$-martingales $(M^n)_{n \geq 1}$ in $\mathcal{M}_0^{c,\mathrm{loc}}(\Omega; E)$ by

$$M^n := I^{W_H}(\xi_{X_n}).$$

Since $\lim_{n \to \infty} X_n = X$ it follows from Lemma 5.4, applied to the differences $X_m - X_n$, that $(M^n)_{n \geq 1}$ is a Cauchy sequence in $L^0(\Omega; C([0,T]; E))$. It follows that $(M^n)_{n \geq 1}$ converges to $M \in L^0(\Omega; C([0,T]; E))$. As a process, $M = (M_t)_{t \in [0,T]}$ is adapted and $M_0 = 0$ almost surely. To show that $M \in \mathcal{M}_0^{c,\mathrm{loc}}(\Omega; E)$ it is now enough to show that $(M_t)_{t \in [0,T]}$ is a local martingale. We claim that

$$M_{\tau_m \wedge t} = M_t^m \qquad \text{almost surely.}$$

This will complete the proof, since it shows that $M$ is a local martingale with localizing sequence $(\tau_m)_{m \geq 1}$. To prove the claim we fix $m \geq 1$. It follows from Proposition 5.3 that for all $n \geq m \geq 1$,

$$\begin{aligned}(5.6)\quad M_{\tau_m \wedge t}^n &= (I^{W_H}(\xi_{X_n}))_{\tau_m \wedge t} = I^{W_H}((\xi_{X_n})_{\tau_m \wedge t}) \\ &= I^{W_H}(\xi_{X_n}(\tau_m \wedge t)) = (I^{W_H}(\xi_{X_m}))_t = M_t^m \qquad \text{almost surely.}\end{aligned}$$

By passing to a subsequence we may assume that $\lim_{n\to\infty} M^n = M$ in $C([0,T];E)$ almost surely. Then also $\lim_{n\to\infty} M^n_{\tau_m \wedge t} = M_{\tau_m \wedge t}$ in $C([0,T];E)$ almost surely, and the claim now follows by letting $n$ tend to infinity in (5.6). It follows that $I^{W_H}(X) := M$ is well defined. At the same time, this argument shows that (5.4) extends to all $X \in L^0_{\mathbb{F}}(\Omega; \gamma(L^2(0,T;H), E))$. This in turn shows that $I^{W_H}$ is continuous.

Next, we extend (5.5) to arbitrary $X \in L^0_{\mathbb{F}}(\Omega; \gamma(L^2(0,T;H), E))$. Let $M = I^{W_H}(\xi_X)$ and define a sequence of stopping times $(\tau_n)_{n \geq 1}$ as
$$\tau_n = \inf\{t \in [0,T] : \|\xi_X(t)\| \geq n\}.$$
By the above results we have, $I^{W_H}(\xi_{X_n}) = M^{\tau_n}$. Applying (5.5) to each $X_n$ and letting $n$ tend to infinity one obtains (5.5) for $X$. From this, we deduce that $I^{W_H}$ has a continuous inverse. This also shows that the mapping $I^{W_H}$ has a closed range in $\mathcal{M}^{c,\text{loc}}_0(\Omega; E)$ and $L^0(\Omega; C([0,T];E))$.

Next assume that $\mathbb{F} = \mathbb{F}^{W_H}$. It suffices to show that the mapping $I^{W_H}$ is surjective. Let $M \in \mathcal{M}^{c,\text{loc}}_0(\Omega, \mathbb{F}^{W_H}; E)$ be arbitrary. We can find a localizing sequence $(\tau_n)_{n \geq 1}$ such that each $M^{\tau_n}$ is a bounded martingale. It follows from the second part of Theorem 3.5 that there is a sequence $(X_n)_{n \geq 1}$ in $L^2_{\mathbb{F}^{W_H}}(\Omega; \gamma(L^2(0,T;H), E))$ such that
$$I^{W_H}(\xi_{X_n}) = M^{\tau_n}.$$
Clearly, $(M^{\tau_n})_{n \geq 1}$ converges to $M$ in $\mathcal{M}^{c,\text{loc}}_0(\Omega, \mathbb{F}^{W_H}; E)$. It follows from Theorem 5.5 that $(X_n)_{n \geq 1}$ is a Cauchy sequence in $L^0_{\mathbb{F}}(\Omega; \gamma(L^2(0,T;H), E))$ and therefore it converges to some $X \in L^0_{\mathbb{F}}(\Omega; \gamma(L^2(0,T;H), E))$. It follows from Theorem 5.5 that $I^{W_H}(X) = M$. □

REMARK 5.6. Proposition 5.3 extends to arbitrary $X \in L^0_{\mathbb{F}}(\Omega; \gamma(L^2(0,T;H), E))$. This may be proved similarly as in Proposition 5.3, but now using Theorem 5.5 for the approximation argument.

The next results on stochastic integration for $H$-valued processes will be used below.

FACTS 5.7. Let $\phi : [0,T] \times \Omega \to H$ be a strongly measurable adapted process such that $\phi \in L^2(0,T;H)$ almost surely. The following results hold:

- The integral process $\int_0^\cdot \phi(t)\,dW_H(t)$ is well defined and belongs to $\mathcal{M}^{c,\text{loc}}_0(\Omega; \mathbb{R})$.
- The quadratic variation process of $\int_0^\cdot \phi(t)\,dW_H(t)$ is given by $\int_0^\cdot \|\phi(t)\|^2\,dt$.
- If $\tau$ is a stopping time, then almost surely for all $t \in [0,T]$ we have
$$\int_0^{\tau \wedge t} \phi(s)\,dW_H(s) = \int_0^t \mathbf{1}_{[0,\tau]}(s)\phi(s)\,dW_H(s).$$

PROPOSITION 5.8. *Let $\Phi:[0,T]\times\Omega\to E$ be an $H$-strongly measurable and adapted process which belongs scalarly to $L^0(\Omega;L^2(0,T;H))$. If there exists a process $\zeta\in L^0(\Omega;C([0,T];E))$ such that for all $x^*\in E^*$ we have*

$$\langle \zeta, x^*\rangle = \int_0^\cdot \Phi^*(t)x^*\,dW_H(t) \qquad \text{in } L^0(\Omega;C([0,T];\mathbb{R})),$$

*then $\zeta$ belongs to $\mathcal{M}_0^{c,\text{loc}}(\Omega;E)$.*

PROOF. Clearly, $\zeta_0 = 0$ almost surely and $\zeta$ is adapted, so it suffices to show $\zeta$ is a local martingale. It is obvious that for all $x^*\in E^*$, $\langle \zeta, x^*\rangle$ is a local martingale. Define a sequence of stopping times $(\tau_n)_{n\geq 1}$ by

$$\tau_n := \inf\{t\in[0,T]: \|\zeta_t\|\geq n\}.$$

By Facts 5.7, for all $x^*\in E^*$ we have

$$\langle \zeta^{\tau_n}, x^*\rangle = \int_0^\cdot \langle \Phi(s), x^*\rangle \mathbf{1}_{[0,\tau_n]}(s)\,dW_H(s) \qquad \text{in } C([0,T];\mathbb{R}) \text{ almost surely}.$$

Since the local martingale on left-hand side is bounded, the Burkholder–Davis–Gundy inequalities and [18], Corollary 17.8, imply that it is a martingale and for all $x^*\in E^*$ and $0\leq s\leq t$ it follows that

$$\langle \mathbb{E}(\zeta_{\tau_n\wedge t}|\mathcal{F}_s), x^*\rangle = \mathbb{E}(\langle \zeta_{\tau_n\wedge t}, x^*\rangle | \mathcal{F}_s) = \langle \zeta_{\tau_n\wedge s}, x^*\rangle$$

almost surely. It follows that for all $0\leq s\leq t$ we have $\mathbb{E}(\zeta_{\tau_n\wedge t}|\mathcal{F}_s) = \zeta_{\tau_n\wedge s}$, so $(\zeta_{\tau_n\wedge t})_{t\in[0,T]}$ is a martingale and $(\zeta_t)_{t\in[0,T]}$ is a local martingale. □

For elementary adapted processes $\Phi:[0,T]\times\Omega\to\mathcal{L}(H,E)$ we define the stochastic integral as an element of $L^0(\Omega;C([0,T];E))$ in the obvious way. The following result extends the integral to a larger class of processes.

THEOREM 5.9. *Let $E$ be a UMD space. For an $H$-strongly measurable and adapted process $\Phi:[0,T]\times\Omega\to\mathcal{L}(H,E)$ which is scalarly in $L^0(\Omega;L^2(0,T;H))$ the following assertions are equivalent:*

(1) *there exists a sequence $(\Phi_n)_{n\geq 1}$ of elementary adapted processes such that:*

   (i) *for all $h\in H$ and $x^*\in E^*$ we have $\lim_{n\to\infty}\langle \Phi_n h, x^*\rangle = \langle \Phi h, x^*\rangle$ in measure on $[0,T]\times\Omega$,*
   (ii) *there exists a process $\zeta\in L^0(\Omega;C([0,T];E))$ such that*

$$\zeta = \lim_{n\to\infty}\int_0^\cdot \Phi_n(t)\,dW_H(t) \qquad \text{in } L^0(\Omega;C([0,T];E));$$

(2) *There exists a process $\zeta \in L^0(\Omega; C([0,T];E))$ such that for all $x^* \in E^*$ we have*

$$\langle \zeta, x^* \rangle = \int_0^{\cdot} \Phi^*(t) x^* \, dW_H(t) \qquad in \ L^0(\Omega; C[0,T]);$$

(3) $\Phi$ *represents an element* $X \in L^0(\Omega; \gamma(L^2(0,T;H), E))$;

(4) *For almost all $\omega \in \Omega$, $\Phi_\omega$ is stochastically integrable with respect to $\tilde{W}_H$.*

*In this situation* $X \in L^0_{\mathbb{F}}(\Omega; \gamma(L^2(0,T;H), E))$ *and*

$$\zeta = I^{W_H}(\xi_X) \qquad in \ L^0(\Omega; C([0,T];E)).$$

A process $\Phi \colon [0,T] \times \Omega \to \mathcal{L}(H, E)$ satisfying the equivalent conditions of the theorem will be called *stochastically integrable* with respect to $W_H$. The process $\zeta = I^{W_H}(\xi_X)$ is called the *stochastic integral process* of $\Phi$ with respect to $W_H$, notation

$$\zeta = \int_0^{\cdot} \Phi(t) \, dW_H(t).$$

It follows from Proposition 5.8 that $\zeta \in \mathcal{M}_0^{\mathrm{c,loc}}(\Omega; E)$.

It is immediate from Proposition 4.3 that if $\Phi \colon [0,T] \times \Omega \to \mathcal{L}(H,E)$ is $L^p$-stochastically integrable for some $p \in (1, \infty)$, then $\Phi$ is stochastically integrable and we have

$$I^{W_H}(\xi_X) = \int_0^{\cdot} \Phi(t) \, dW_H(t),$$

where $X \in L^p(\Omega; \gamma(L^2(0,T;H), E))$ is represented by $\Phi$.

REMARK 5.10. Under the assumptions as stated, condition (1) is equivalent to:

(1)′ There exists a sequence $(\Phi_n)_{n \geq 1}$ of elementary adapted processes such that:

(i) for all $h \in H$ we have $\lim_{n \to \infty} \Phi_n h = \Phi h$ in measure on $[0,T] \times \Omega$;
(ii) there exists an $\eta \in L^0(\Omega; C([0,T];E))$ such that

$$\eta = \lim_{n \to \infty} \int_0^{\cdot} \Phi_n(t) \, dW_H(t) \qquad in \ L^0(\Omega; C([0,T];E)).$$

PROOF OF THEOREM 5.9. First note that (i) and (ii) of part (1), combined with [18], Proposition 17.6, imply that in (i) we have convergence in $L^0(\Omega; L^2(0,T;H))$.

$(1) \Rightarrow (3)$: Let $\Phi_n$ represent $X_n \in L^0(\Omega; \gamma(L^2(0,T;H), E))$. By (ii) and Lemma 5.4, these elements define a Cauchy sequence in $L^0(\Omega; \gamma(L^2(0,T;H),$ $E$)). Let $X \in L^0(\Omega; \gamma(L^2(0,T;H), E))$ be the limit. Since each $X_n$ is elementary adapted we have $X \in L^0_{\mathbb{F}}(\Omega; \gamma(L^2(0,T;H), E))$, and with property (i) it follows that

$$\langle X, x^* \rangle = \lim_{n \to \infty} \langle X_n, x^* \rangle = \lim_{n \to \infty} \Phi_n^* x^* = \Phi^* x^* \qquad \text{in } L^0(\Omega; L^2(0,T;H)).$$

Hence, $\Phi$ represents $X$.

(3) $\Rightarrow$ (4): It follows from Lemma 2.7 that for almost all $\omega \in \Omega$, $\Phi_\omega$ is represented by $X(\omega)$. The result now follows from Proposition 3.2.

(4) $\Rightarrow$ (3): Let $N$ be a null set such that $\Phi_\omega$ is stochastically integrable with respect to $\tilde{W}_H$ for all $\omega \in \complement N$. Proposition 3.2 assures that for such $\omega$ we may define $X(\omega) \in \gamma(L^2(0,T;H), E)$ defined by

$$X(\omega)f = \int_0^T \Phi(t, \omega) f(t) \, dt.$$

An application of Remark 2.8 shows that the resulting random variable $X : \Omega \to \gamma(L^2(0,T;H), E)$ is strongly measurable. This proves (2).

(3) $\Rightarrow$ (1): This may be proved in the same way as Theorem 3.6, this time using Theorem 5.5.

(1) $\Rightarrow$ (2): This is clear.

(2) $\Rightarrow$ (1): It follows from Proposition 5.8 that $\zeta \in \mathcal{M}_0^{c,\text{loc}}(\Omega; E)$. Let $(\tau_n)_{n \geq 1}$ be a localizing sequence such that each $\zeta^{\tau_n}$ is bounded. It follows from the assumptions and Facts 5.7 that for all $n$ and all $x^* \in E^*$ we have

$$\langle \zeta^{\tau_n}, x^* \rangle = \int_0^\cdot \mathbf{1}_{[0, \tau_n]}(t) \Phi^*(t) x^* \, dW_H(t) \qquad \text{almost surely.}$$

By the Burkholder–Davis–Gundy inequalities, each $\mathbf{1}_{[0,\tau_n]}\Phi$ is scalarly in $L^2(\Omega; L^2(0,T;H))$. In particular,

$$\langle \zeta_{\tau_n}, x^* \rangle = \int_0^T \mathbf{1}_{[0,\tau_n]}(t) \Phi^*(t) x^* \, dW_H(t) \qquad \text{in } L^2(\Omega).$$

By Theorem 3.6, each $\mathbf{1}_{[0,\tau_n]}\Phi$ is $L^2$-stochastically integrable with integral $\zeta_{\tau_n}$. With Theorem 3.6 we find elementary adapted processes $(\Phi_n)_{n \geq 1}$ such that

$$\left\| \zeta_{\tau_n} - \int_0^T \Phi_n(t) \, dW_H(t) \right\|_{L^2(\Omega;E)} < \frac{1}{n}.$$

Doob's maximal inequality implies that

$$\left\| \zeta^{\tau_n} - \int_0^\cdot \Phi_n(t) \, dW_H(t) \right\|_{L^2(\Omega; C([0,T];E))} \leq \frac{2}{n}.$$

It follows that

$$\left\|\zeta - \int_0^\cdot \Phi_n(t)\,dW_H(t)\right\|_{L^0(\Omega;C([0,T];E))}$$
$$\leq \|\zeta - \zeta^{\tau_n}\|_{L^0(\Omega;C([0,T];E))} + \left\|\zeta^{\tau_n} - \int_0^\cdot \Phi_n(t)\,dW_H(t)\right\|_{L^0(\Omega;C([0,T];E))}$$
$$\leq \|\zeta - \zeta^{\tau_n}\|_{L^0(\Omega;C([0,T];E))} + \frac{2}{n}.$$

The latter clearly converges to 0 as $n$ tends to infinity. This gives (ii). Now choose $x^* \in E^*$ arbitrary. In view of

$$\int_0^\cdot \Phi^*(t)x^*\,dW_H(t) = \lim_{n\to\infty} \int_0^\cdot \Phi_n^*(t)x^*\,dW_H(t) \qquad \text{in } L^0(\Omega;C([0,T]))$$

from [18], Proposition 17.6, we obtain (i). □

REMARK 5.11. As was the case in Remark 3.8, if the filtration $\mathbb{F}$ is assumed to be the augmented Brownian filtration $\mathcal{F}_T^{W_H}$, then the equivalence $(1) \Leftrightarrow (2)$ is true for every real Banach space $E$. This may be proved by a stopping time argument as in the proof of $(2) \Rightarrow (1)$.

Our next objective is a generalization Theorem 4.4.

THEOREM 5.12 (Burkholder–Davis–Gundy inequalities). *Let $E$ be a UMD space and fix $p \in (1,\infty)$. If $\Phi:[0,T] \times \Omega \to \mathcal{L}(H,E)$ is $H$-strongly measurable and adapted and stochastically integrable, then*

$$\mathbb{E} \sup_{t\in[0,T]} \left\|\int_0^t \Phi(s)\,dW_H(s)\right\|^p \eqsim_{p,E} \mathbb{E}\|X\|_{\gamma(L^2(0,T;H),E)}^p,$$

*where $X \in L_\mathbb{F}^0(\Omega;\gamma(L^2(0,T;H),E))$ is the element represented by $\Phi$.*

This is understood in the sense that the left-hand side is finite if and only if the right-hand side is finite, in which case the estimates hold with constants only depending on $p$ and $E$.

PROOF OF THEOREM 5.12. First assume that the left-hand side is finite. Define a sequence of stopping times $(\tau_n)_{n\geq 1}$ by

$$\tau_n = \inf\{t \in [0,T] : \|\xi_X(t)\|_{\gamma(L^2(0,T;H),E)} \geq n\}.$$

Observe that $\xi_X(\tau_n) \in L_\mathbb{F}^p(\Omega;\gamma(L^2(0,T;H),E))$ and that it is represented by $\Phi \mathbf{1}_{[0,\tau_n]}$. From Theorem 3.6 we deduce that $\Phi \mathbf{1}_{[0,\tau_n]}$ is $L^p$-stochastically integrable. Combining the identity

$$\int_0^{\tau_n} \Phi(t)\,dW_H(t) = \int_0^T \mathbf{1}_{[0,\tau_n]}(t)\Phi(t)\,dW_H(t)$$

which follows for instance from Theorem 5.9(1), with the dominated convergence theorem (here we use the assumption) and Fatou's lemma, we obtain

$$\mathbb{E}\left\|\int_0^T \Phi(t)\,dW_H(t)\right\|^p = \lim_{n\to\infty}\mathbb{E}\left\|\int_0^T \mathbf{1}_{[0,\tau_n]}(t)\Phi(t)\,dW_H(t)\right\|^p$$

$$\eqsim_{p,E} \liminf \|\xi_X(\tau_n)\|^p_{L^p(\Omega;\gamma(L^2(0,T;H),E))}$$

$$\geq \|X\|^p_{L^p(\Omega;\gamma(L^2(0,T;H),E))}.$$

This shows that $X \in L^p(\Omega; \gamma(L^2(0,T;H), E))$, and by Theorem 3.6 that $\Phi$ is $L^p$-stochastically integrable. The result now follows from Theorem 4.4.

If the right-hand side is finite, then $\Phi$ is $L^p$-stochastically integrable by Theorem 3.6 and therefore the left-hand side is finite by Theorem 4.4. $\square$

In the real-valued case, a similar estimates holds for all $0 < p < \infty$. We do not know whether Theorem 5.12 extends to all $0 < p < \infty$ (or even just to $p = 1$).

We have the following extension of Itô's representation theorem for Brownian martingales to UMD Banach spaces.

THEOREM 5.13 (Representation of UMD-valued Brownian local martingales). *Let $E$ be a UMD space. Then every $E$-valued local martingale $M := (M_t)_{t\in[0,T]}$ adapted to the augmented filtration $\mathbb{F}^{W_H}$ has a continuous version and there exists a unique $X \in L^0_{\mathbb{F}}(\Omega; \gamma(L^2(0,T;H), E))$ such that*

$$M = M_0 + I^{W_H}(\xi_X).$$

PROOF. We may assume $M_0 = 0$. By Theorem 5.5 it suffices to show that $M$ has a continuous version. This can be seen in the same way as in the real case (cf. [18], Theorem 18.10). $\square$

For UMD spaces $E$ with cotype 2 recall that $\gamma(L^2(0,T;H), E) \hookrightarrow L^2(0,T; \gamma(H,E))$. Hence every $X \in L^0(\Omega; \gamma(L^2(0,T), E))$ can be represented by a process $\Phi \in L^0(\Omega; L^2(0,T; \gamma(H, E)))$. In this case, the above representation takes the form

$$M = M_0 + \int_0^{(\cdot)} \Phi(t)\,dW_H(t).$$

For $M$-type 2 Banach spaces $E$, a representation theorem for martingales as stochastic integrals with respect to $H$-cylindrical Brownian motions can be found in [31], Chapter 2.

**Acknowledgment.** M. Veraar gratefully acknowledges support from the Marie Curie Fellowship Program for a stay at the TU Karlsruhe.

# REFERENCES

[1] BELOPOLSKAYA, YA. I. and DALETSKIĬ, YU. L. (1990). *Stochastic Equations and Differential Geometry.* Kluwer Academic, Dordrecht. MR1050097

[2] BOURGAIN, J. (1983). Some remarks on Banach spaces in which martingale difference sequences are unconditional. *Ark. Mat.* **21** 163–168. MR0727340

[3] BOURGAIN, J. (1986). Vector-valued singular integrals and the $H^1$-BMO duality. In *Probability Theory and Harmonic Analysis* (*Cleveland, Ohio, 1983*) 1–19. *Monogr. Textbooks Pure Appl. Math.* **98**. Dekker, New York. MR0830227

[4] BRZEŹNIAK, Z. (1995). Stochastic partial differential equations in M-type 2 Banach spaces. *Potential Anal.* **4** 1–45. MR1313905

[5] BRZEŹNIAK, Z. (1997). On stochastic convolutions in Banach spaces and applications. *Stoch. Stoch. Rep.* **61** 245–295. MR1488138

[6] BRZEŹNIAK, Z. (2003). Some remarks on stochastic integration in 2-smooth Banach spaces. In *Probabilistic Methods in Fluids* (I. M. Davies, A. Truman et al., eds.) 48–69. World Scientific, New Jersey. MR2083364

[7] BRZEŹNIAK, Z. and VAN NEERVEN, J. M. A. M. (2000). Stochastic convolution in separable Banach spaces and the stochastic linear Cauchy problem. *Studia Math.* **143** 43–74. MR1814480

[8] BURKHOLDER, D. L. (2001). Martingales and singular integrals in Banach spaces. In *Handbook of the Geometry of Banach Spaces* 233–269. North-Holland, Amsterdam. MR1863694

[9] CHUNG, K. L. and WILLIAMS, R. J. (1990). *Introduction to Stochastic Integration*, 2nd ed. Birkhäuser, Boston. MR1102676

[10] CLÉMENT, PH. P. J. E., DE PAGTER, B., SUKOCHEV, F. A. and WITVLIET, H. (2000). Schauder decompositions and multiplier theorems. *Studia Math.* **138** 135–163. MR1749077

[11] DETTWEILER, E. (1989). On the martingale problem for Banach space valued stochastic differential equations. *J. Theoret. Probab.* **2** 159–191. MR0987575

[12] DETTWEILER, E. (1991). Stochastic integration relative to Brownian motion on a general Banach space. *Doğa Mat.* **15** 58–97. MR1115509

[13] DIESTEL, J., JARCHOW, H. and TONGE, A. (1995). *Absolutely Summing Operators.* Cambridge Univ. Press. MR1342297

[14] EDGAR, G. A. (1977). Measurability in a Banach space. *Indiana Univ. Math. J.* **26** 663–677. MR0487448

[15] GARLING, D. J. H. (1986). Brownian motion and UMD-spaces. *Probability and Banach Spaces* (*Zaragoza, 1985*). *Lecture Notes in Math.* **1221** 36–49. Springer, Berlin. MR0875006

[16] GARLING, D. J. H. (1990). Random martingale transform inequalities. In *Probability in Banach Spaces VI* (*Sandbjerg, 1986*) 101–119. Birkhäuser, Boston. MR1056706

[17] HITCZENKO, P. (1988). On tangent sequences of UMD-space valued random vectors. Unpublished manuscript. Warsaw.

[18] KALLENBERG, O. (2002). *Foundations of Modern Probability*, 2nd ed. Springer, New York. MR1876169

[19] KALTON, N. J. and WEIS, L. The $H^\infty$-functional calculus and square function estimates. Preprint.

[20] KARATZAS, I. and SHREVE, S. E. (1991). *Brownian Motion and Stochastic Calculus*, 2nd ed. Springer, New York. MR1121940

[21] KUNITA, H. (1970). Stochastic integrals based on martingales taking values in Hilbert space. *Nagoya Math. J.* **38** 41–52. MR0264754

[22] LEDOUX, M. and TALAGRAND, M. (1991). *Probability in Banach Spaces*. Springer, Berlin. [MR1102015](MR1102015)

[23] MAMPORIA, B. (2004). On the existence and uniqueness of a solution to a stochastic differential equation in a Banach space. *Georgian Math. J.* **11** 515–526. [MR2081742](MR2081742)

[24] MCCONNELL, T. R. (1989). Decoupling and stochastic integration in UMD Banach spaces. *Probab. Math. Statist.* **10** 283–295. [MR1057936](MR1057936)

[25] MONTGOMERY-SMITH, S. (1998). Concrete representation of martingales. *Electron. J. Probab.* **3** 15. [MR1658686](MR1658686)

[26] VAN NEERVEN, J. M. A. M., VERAAR, M. C. and WEIS, L. (2006). Itô's formula in UMD Banach spaces and regularity of solutions of the Zakai equation. Unpublished manuscript.

[27] VAN NEERVEN, J. M. A. M., VERAAR, M. C. and WEIS, L. (2006). Stochastic evolution equations in UMD Banach spaces. Unpublished manuscript.

[28] VAN NEERVEN, J. M. A. M. and WEIS, L. (2005). Stochastic integration of functions with values in a Banach space. *Studia Math.* **166** 131–170. [MR2109586](MR2109586)

[29] VAN NEERVEN, J. M. A. M. and WEIS, L. (2005). Weak limits and integrals of Gaussian covariances in Banach spaces. *Probab. Math. Statist.* **25** 55–74. [MR2211356](MR2211356)

[30] NEIDHARDT, A. L. (1978). Stochastic integrals in 2-uniformly smooth Banach spaces. Ph.D. dissertation, Univ. Wisconsin.

[31] ONDREJÁT, M. (2003). Équations d'Évolution Stochastiques dans les Espaces de Banach. Ph.D. thesis. Inst. Élie Cartan, Nancy, and Charles Univ., Prague.

[32] PISIER, G. (1975). Martingales with values in uniformly convex spaces. *Israel J. Math.* **20** 326–350. [MR0394135](MR0394135)

[33] PISIER, G. (1986). Probabilistic methods in the geometry of Banach spaces. *Probability and Analysis* (*Varenna, 1985*). *Lecture Notes in Math.* **1206** 167–241. Springer, Berlin. [MR0864714](MR0864714)

[34] RUBIO DE FRANCIA, J. L. (1986). Martingale and integral transforms of Banach space valued functions. *Probability and Banach Spaces* (*Zaragoza, 1985*). *Lecture Notes in Math.* **1221** 195–222. Springer, Berlin. [MR0875011](MR0875011)

[35] ROSIŃSKI, J. (1987). Bilinear random integrals. *Dissertationes Math.* **259** 71. [MR0888463](MR0888463)

[36] ROSIŃSKI, J. and SUCHANECKI, Z. (1980). On the space of vector-valued functions integrable with respect to the white noise. *Colloq. Math.* **43** 183–201. [MR0615985](MR0615985)

[37] VAKHANIA, N. N., TARIELADZE, V. I. and CHOBANYAN, S. A. (1987). *Probability Distributions in Banach Spaces*. Reidel, Dordrecht. [MR1435288](MR1435288)

J. M. A. M. VAN NEERVEN
M. C. VERAAR
DELFT INSTITUTE OF
  APPLIED MATHEMATICS
DELFT UNIVERSITY OF TECHNOLOGY
P.O. BOX 5031
2600 GA DELFT
THE NETHERLANDS
E-MAIL: [J.M.A.M.vanNeerven@tudelft.nl](mailto:J.M.A.M.vanNeerven@tudelft.nl)
        [M.C.Veraar@tudelft.nl](mailto:M.C.Veraar@tudelft.nl)
URL: [http://fa.its.tudelft.nl/~neerven/](http://fa.its.tudelft.nl/~neerven/)
     [http://fa.its.tudelft.nl/~veraar/](http://fa.its.tudelft.nl/~veraar/)

L. WEIS
MATHEMATISCHES INSTITUT I
TECHNISCHE UNIVERSITÄT KARLSRUHE
D-76128 KARLSRUHE
GERMANY
E-MAIL: [Lutz.Weis@math.uni-karlsruhe.de](mailto:Lutz.Weis@math.uni-karlsruhe.de)
URL: [http://www.mathematik.uni-karlsruhe.de/mi1weis/~weis/en](http://www.mathematik.uni-karlsruhe.de/mi1weis/~weis/en)